\documentclass[10pt]{article}

\usepackage{amssymb}
\usepackage{enumerate}
\usepackage{amsmath}
\newtheorem{theorem}{Theorem}[section]
\newtheorem{corollary}[theorem]{Corollary}
\newtheorem{lemma}[theorem]{Lemma}

\newtheorem{property}[theorem]{Property}
\newtheorem{definition}[theorem]{Definition}
\newtheorem{remark}[theorem]{Remark}
\begin{document}

%----------------------------------------------------------------Title------------------------------------------------------------------------------------%

\title{\bf Coherent pairs of measures and Markov-Bernstein inequalities}%
\maketitle

%-------------------------------------------------------------Authors+Adress------------------------------------------------------------------------------%

\centerline{\bf Andr\'e Draux }
\centerline{\it Normandie Univ, INSA Rouen, LMI, 76000 Rouen, France}
%\centerline{\it Laboratoire de Math\'ematique de l'INSA de Rouen - LMI (EA 3226 - FR CNRS 3335),}
\centerline{\it Campus de Saint-\'Etienne-du-Rouvray,}
\centerline{\it 685 Avenue de l'universit\'e, BP 8}
\centerline{\it F-76801 Saint-\'Etienne-du-Rouvray Cedex, France}
\centerline{\it e-mail: andre.draux@insa-rouen.fr}
\bigskip

%\author{Andr\'e Draux}%
%\address{Normandie Universit\'e\\Laboratoire de Math\'ematique de l'INSA de Rouen - LMI (EA 3226 - FR CNRS 3335), Campus de Saint-\'Etienne-du-Rouvray, 685 Avenue de l'universit\'e, BP 8, F-76801 Saint-\'Etienne-du-Rouvray Cedex, France}
%\email{andre.draux@insa-rouen.fr}%
% ----------------------------------------------------------------
% AMS-LaTeX Paper ************************************************
% **** -----------------------------------------------------------
%\documentclass[10pt]{article}
%\documentclass{amsart}
% --------------------------------------------------------Abstract--------------------------------------------------------------------

%
 
\maketitle
\begin{abstract}
All the coherent pairs of measures associated to linear functionals $c_0$ and $c_1$, introduced by Iserles et al in 1991, have been given by Meijer in 1997. There exist seven kinds of coherent pairs. All these cases are explored in order to give three term recurrence relations satisfied by polynomials. The smallest zero $\mu_{1,n}$ of each of them of degree $n$ has a link with the Markov-Bernstein constant $M_n$ appearing in the following Markov-Bernstein inequalities:
$$
c_1((p^\prime)^2) \le M_n^2 c_0(p^2), \quad \forall p \in \mathcal{P}_n,
$$
where $M_n=\frac{1}{\sqrt{\mu_{1,n}}}$.\\
The seven kinds of three term recurrence relations are given. In the case where $c_0 =e^{-x} dx+\delta(0)$ and $c_1 =e^{-x} dx$, explicit upper and lower bounds are given for $\mu_{1,n}$, and the asymptotic behavior of the corresponding Markov-Bernstein constant is stated. Except in a part of one case, $\lim_{n \to \infty} \mu_{1,n}=0$ is proved in all the cases.
\end{abstract}
%-------------------------------------------------------Keywords--------------------------------------------------------------------------------------------%

\medskip


%**************************************************************************************************************************

%_______________________________________________section 1 ___________________________

%**************************************************************************************************************************

\section{Introduction}

The Markov-Bernstein inequalities in $L^2$ norm that imply two measures, have the following form
$$
c_1((p^\prime)^2) \le M_n^2 c_0(p^2), \quad \forall p \in \mathcal{P}_n,
$$
where $c_0$ and $c_1$ denote the two linear functionals associated to the both measures. $\mathcal{P}_n$ is the vector space of real polynomials in one variable of degree at most $n$. $M_n$ is called Markov-Bernstein constant. These inequalities are always related to an eigenvalue problem of a positive definite symmetric matrix (see \cite{MMR}, \cite{DEH2} for a general presentation). For any measures this matrix generally is full. But for classical measures (Hermite, Laguerre-Sonin, Jacobi) this matrix has a particular form. It is diagonal for all these measures (see \cite{DEH1}) if $c_0=c_1= \int_{-\infty}^{+\infty}. e^{-x^2}dx$ in the Hermite case, if $c_0=c^\alpha= \int_{0}^{+\infty}. x^\alpha e^{-x}dx$ with $\alpha >-1$ and $c_1=c^{\alpha+1}$ in the Laguerre-Sonin case, if $c_0=c^{(\alpha,\beta)}= \int_{-1}^{+1}.(1- x)^\alpha (1+ x)^\beta dx$ with $\alpha >-1$ and $\beta >-1$, and $c_1=c^{(\alpha+1,\beta+1)}$ in the Jacobi case. It is tridiagonal (see \cite{MMR}, \cite{DEH2}) if $c_0=c_1=c^\alpha$ in the Laguerre-Sonin case. If $c_0=c_1=c^{(\alpha,\alpha)}$ that is the Gegenbauer case, two tridiagonal matrices are obtained (see \cite{MMR}, \cite{DEH2}). It is a five-diagonal matrix if $c_0=c_1=c^{(\alpha,\beta)}$ in the Jacobi case (see \cite{DEH2}, \cite{DMS}). In the discrete case, if $c_0$ and $c_1$ are identical and correspond to the Meixner measure, the matrix is also tridiagonal (see \cite{DRA1}, \cite{ADT}).
When the matrix is tridiagonal, its characteristic polynomial satisfy a three term recurrence relation. To find the Markov-Bernstein constant is equivalent to find the smallest zero $\mu_{1,n}$ of the polynomial of degree $n$ satisfying this three term recurrence relation, and $M_n = 1/\sqrt{\mu_{1,n}}$. The characteristic polynomial in the Jacobi case  $c_0=c_1=c^{(\alpha,\beta)}$ provides a six-term recurrence relation (see \cite{DMS}). Thus, to obtain a tridiagonal matrix is exceptional. There exists another case for which a tridiagonal matrix is obtained: it is the one of coherent pairs $(c_0,c_1)$.

The notion of coherent pairs of measures was introduced for first time by Iserles et al \cite{IKNS} in 1991. These measures are defined as follows

%definition
\begin{definition}\label{D1}
Let $c_0$ and $c_1$ be two quasi definite linear functionals. Let $\{P_n\}_{n\ge 0}$ (resp. $\{T_n\}_{n\ge 0}$) be the sequence of monic orthogonal polynomials with respect to  $c_0$ (resp. $c_1$). $(c_0,c_1)$ is called coherent pair if and only if there exists a sequence $\{\sigma_n\}_{n\ge 1}$, $\sigma_n \in \mathbb{R}$, $\sigma_n \ne 0$, such that
\begin{equation}\label{eq:1.1}
T_n(x) = \frac{P^\prime_{n+1} (x)}{n+1}-\sigma_n \frac{P^\prime_{n} (x)}{n}, \quad \forall n \ge 1.
\end{equation}
\end{definition}
All the kinds of coherent pairs $(c_0,c_1)$ were described by Meijer \cite{MEI} in 1997. They correspond to the relations (3.12) to (3.18) of his paper. One of both functionals is classical (Laguerre-Sonin or Jacobi). The seven cases contained in \cite{MEI} are given below.
\begin{enumerate}
\item {\bf Laguerre case:} $\Omega=]0,+\infty[$
   \begin{enumerate}
   \item $c_0$ corresponds to the measure $(x-\xi)x^{\alpha -1} e^{-x}dx$ with $\alpha >0$ and $\xi < 0$.\\
    $c_1$ corresponds to the Laguerre-Sonin measure $x^{\alpha } e^{-x}dx$.
    \item $c_0$ corresponds to the measure $e^{-x}dx+M\delta(0)$ with $M \ge0$. $\delta$ is the Dirac measure.\\
    $c_1$ corresponds to the Laguerre measure $e^{-x}dx$.
    \item $c_0$ corresponds to the Laguerre-Sonin measure $x^\alpha e^{-x}dx$ with $\alpha >-1$.\\
    $c_1$ corresponds to the measure $\frac{x^{\alpha +1}}{x-\xi}e^{-x}dx +M \delta(\xi)$ with $\xi \le 0$ and $M\ge 0$.
   \end{enumerate}
\item {\bf Jacobi case:} $\Omega=]-1,+1[$
   \begin{enumerate}
   \item $c_0$ corresponds to the measure  $\mid x-\xi \mid (1-x)^{\alpha-1}(1+x)^{\beta-1}dx$ with $\alpha >0$, $\beta >0$ and $\mid \xi \mid > 1$.\\
    $c_1$ corresponds to the Jacobi measure $(1-x)^{\alpha}(1+x)^{\beta}dx$.
    \item $c_0$ corresponds to the measure $(1+x)^{\beta-1}dx+M \delta(1)$ with $\beta >0$ and $M \ge 0$.\\
    $c_1$ corresponds to the Jacobi measure $(1+x)^{\beta}dx$.
    \item $c_0$ corresponds to the measure  $(1-x)^{\alpha-1}  dx +M \delta(-1)$ with $\alpha >0$ and $M \ge 0$.\\
    $c_1$ corresponds to the Jacobi measure $(1-x)^{\alpha}dx$.
    \item $c_0$ corresponds to the Jacobi measure   $(1-x)^{\alpha}(1+x)^{\beta}dx$ with $\alpha >-1$ and $\beta >-1$.\\
     $c_1$ corresponds to the measure $ \frac{1}{\mid x-\xi \mid}(1-x)^{\alpha+1}(1+x)^{\beta+1}dx+M \delta(\xi)$ with $\mid \xi \mid \ge 1$ and $M \ge 0$.
   \end{enumerate}
\end{enumerate}
In every previous cases one measure depends on one parameter $\xi$ or $M$, or both $\xi$ and $M$. Our aim is to fix these parameters, and therefore to fix the different measures. From these measures we want to give, in an explicit form only depending on the fixed parameters, the two sequences of polynomials $\{P_n\}_{n\ge 0}$ and $\{T_n\}_{n\ge 0}$ as well as their $L^2$ norms, and the sequence $\{\sigma_n\}_{n\ge 1}$. After that part we want to obtain, also in an explicit form only depending on the fixed parameters, the three term recurrence relation linked with the study of the positivity of the bilinear functional $a_\lambda(p,q) = c_0(pq)+ \lambda c_1(p^\prime q^\prime)$ for $p,q \in \mathcal{P}$, where $\mathcal{P}$ is the vector space of real polynomials in one variable.\\
We begin to give the general expression of such a three term recurrence relation.\\
Let us denote the square norm of $P_n$ (resp. $T_n$) by $k_n^{(0)}$ (resp. $k_n^{(1)}$).
$$
k_n^{(0)}= c_0(P_n^2) \mbox{ and } k_n^{(1)}= c_1(T_n^2), \quad \forall n \ge 0.
$$
Like in several papers (see \cite{DEH1}, \cite{DEH2}, \cite{DRA}) a particular basis of monic polynomials $\{R_i\}_{i\ge 0}^n$ of $\mathcal{P}_n$ is defined from the following relations: 
\begin{eqnarray*}
R^\prime_0(x) &= &0,\\
R^\prime_1(x) &= & P^\prime_1(x) =1,\\
R^\prime_i(x) &= &i T_{i-1}(x) = P^\prime_i(x) -\sigma_{i-1} \frac{i}{i-1}P^\prime_{i-1}(x), \quad  i = 2,\ldots,n.
\end{eqnarray*}
From these relations the polynomials $R_i$ are chosen as follows:
\begin{eqnarray*}
R_0(x)&= &1,\\
R_1(x) &= & P_1(x),\\
R_i (x)&= &P_i(x) -\sigma_{i-1} \frac{i}{i-1}P_{i-1}(x), \quad  i = 2,\ldots,n.\\
\end{eqnarray*}
Let $p$ be a polynomial belonging to $\mathcal{P}_n$, written in this particular basis: $p(x)=Ê\sum_{i=0}^n y_i R_i(x)$. $y_i \in \mathbb{R}$, $i=0,\ldots,n$. Let $y$ be the vector of $\mathbb{R}^{n+1}$ of components  $y_i$, $i=0,\ldots,n$. Then
$$
a_\lambda(p,p) = y^T K_{n,0}^{(0)} y + \lambda  y^T K_{n,0}^{(1)} y
$$
where $ K_{n,0}^{(0)}$ is the $(n+1) \times (n+1)$ positive definite symmetric tridiagonal matrix
$$
K_{n,0}^{(0)} = (c_0(R_j R_i))_{i,j=0}^n
$$
and $K_{n,0}^{(1)}$ is the $(n+1) \times (n+1)$ positive semidefinite diagonal matrix
$$
K_{n,0}^{(1)}= (c_1(R^\prime_j R^\prime_i))_{i,j=0}^n =(ijc_1(T_{j-1} T_{i-1}))_{i,j=0}^n.
$$
By convention $T_{-1}=0$.\\
Hence the content of both matrices $K_{n,0}^{(0)} $ and $K_{n,0}^{(1)} $.
\begin{eqnarray*}
K_{n,0}^{(1)}(i,i) & = & i^2 c_1(T_{i-1}^2) = i^2 k_{i-1}^{(1)}, \quad i=0,\ldots,n.\\
K_{n,0}^{(0)}(i,i)&=&k_{i}^{(0)}, \quad i=0,1,\\
&=& k_{i}^{(0)} + \sigma_{i-1}^2 (\frac{i}{i-1})^2 k_{i-1}^{(0)}, \quad   i = 2,\ldots,n.\\
K_{n,0}^{(0)}(i+1,i)=K_{n,0}^{(0)}(i,i+1)&=& 0 \quad \mbox{ if } i=0,\\
&=&-\sigma_{i} (\frac{i+1}{i}) k_{i}^{(0)}, \quad   i = 1,\ldots,n-1.
\end{eqnarray*}
We have the following obvious property
%property
\begin{property}\label{P2}
$a_\lambda(p,p) >0$ $\forall p \in \mathcal{P}_n$ if and only if $K_{n,0}^{(0)}  + \lambda   K_{n,0}^{(1)}$  is a positive definite symmetric matrix.
\end{property}
The first row and the first column of $K_{n,0}^{(0)}  + \lambda   K_{n,0}^{(1)}$ only  contain one entry in row 0 and column 0, equal to $k_{0}^{(0)}>0$. We may reduce the size of $K_{n,0}^{(0)}  + \lambda   K_{n,0}^{(1)}$ and only keep 
\begin{eqnarray*}
K_{n,1}^{(0)} &=& (c_0(R_j R_i))_{i,j=1}^n,\\
K_{n,1}^{(1)} &=& (c_1(R^\prime_j R^\prime_i))_{i,j=1}^n.
\end{eqnarray*}
Then, Property \ref{P2} is equivalent to  Property \ref{P3}.
%property
\begin{property}\label{P3}
$a_\lambda(p,p) >0$ $\forall p \in \mathcal{P}_n$ if and only if $K_{n,1}^{(0)}  + \lambda   K_{n,1}^{(1)}$  is a positive definite symmetric matrix.
\end{property}
A transformation of the matrix $K_{n,1}^{(0)}  + \lambda   K_{n,1}^{(1)}$ uses the Cholesky decomposition of $K_{n,1}^{(0)}$.
$$
K_{n,1}^{(0)} = G_{n,1}^{(0)} (G_{n,1}^{(0)})^T
$$
where
\begin{eqnarray*}
G_{n,1}^{(0)}(i,i)&=&\sqrt{k_{i}^{(0)}} \quad \mbox{ if } i\ge1,\\
G_{n,1}^{(0)}(i,i-1)&=& -\sigma_{i-1} \frac{i}{i-1} \sqrt{k_{i-1}^{(0)}},  \quad   i = 2,\ldots,n.
\end{eqnarray*}
Then, we have
$$
K_{n,1}^{(0)} + \lambda K_{n,1}^{(1)} = (K_{n,1}^{(1)})^{1/2}\left((K_{n,1}^{(1)})^{-1/2}G_{n,1}^{(0)} (G_{n,1}^{(0)})^T(K_{n,1}^{(1)})^{-1/2} +\lambda I\right)(K_{n,1}^{(1)})^{1/2}.
$$
Let us denote $(K_{n,1}^{(1)})^{-1/2}G_{n,1}^{(0)} (G_{n,1}^{(0)})^T(K_{n,1}^{(1)})^{-1/2}$ by $\tilde K_{n,1}$ which is a $n \times n$ positive definite symmetric tridiagonal matrix. Its entries are
\begin{eqnarray*}
\tilde K_{n,1}(i,i)&=&\frac{k_{1}^{(0)}}{k_{0}^{(1)}} \quad \mbox{ if } i=1,\\
&=& \frac{k_{i}^{(0)}}{i^2 k_{i-1}^{(1)}}+ \frac{\sigma_{i-1}^2 }{({i-1})^2}\frac{ k_{i-1}^{(0)}}{k_{i-1}^{(1)}}, \quad   i = 2,\ldots,n.\\
\tilde K_{n,1}(i+1,i)=\tilde K_{n,1}(i,i+1) &=& -\frac{\sigma_i k_{i}^{(0)}}{i^2 \sqrt{k_{i-1}^{(1)}k_{i}^{(1)}}}, \quad   i = 1,\ldots,n-1.
\end{eqnarray*}
Property \ref{P3} implies that $\tilde K_{n,1} + \lambda I$ is positive definite. Thus $-\lambda$ has to be smaller than the smallest eigenvalue of $\tilde K_{n,1}$. The eigenvalues of $\tilde K_{n,1}$ also are the zeros of the characteristic polynomial $A_n(\lambda)=\det(\lambda I - \tilde K_{n,1})$. Of course the eigenvalues of $\tilde K_{n,1}$ are positive, so are the zeros of $A_n(\lambda)$. Since $\tilde K_{n,1}$ is a tridiagonal matrix, the polynomials $A_n(\lambda)$ satisfy a three term recurrence relation.
\begin{equation}\label{eq:6.1}
A_n (\lambda) = (\lambda - \frac{k_{n}^{(0)}}{n^2 k_{n-1}^{(1)}}- \frac{\sigma_{n-1}^2 }{({n-1})^2}\frac{ k_{n-1}^{(0)}}{k_{n-1}^{(1)}}) A_{n-1} (\lambda)- \frac{\sigma_{n-1}^2 (k_{n-1}^{(0)})^2}{(n-1)^4 k_{n-1}^{(1)}k_{n-2}^{(1)}} A_{n-2} (\lambda)
\end{equation}
with $A_{0} (\lambda) =1$ and $A_{1} (\lambda) = \lambda - \frac{k_{1}^{(0)}}{k_{0}^{(1)}}$.\\
This three term recurrence relation was also given by P\'erez \cite{P} in her thesis (Relation (5.4.13) page 157), but obtained by another way.\\
The smallest zero $\mu_{1,n}$ of $A_n (\lambda)$ provides the following Markov-Bernstein inequality
%theorem
\begin{theorem}\label{1.4}
The following Markov-Bernstein inequality holds
\begin{equation}\label{eq:6.1b}
c_1((p^\prime)^2) \le \frac{1}{\mu_{1,n}} c_0(p^2), \quad \forall p \in \mathcal{P}_n
\end{equation}
where $\mu_{1,n}$ is the smallest zero $\mu_{1,n}$ of the polynomial $A_n (\lambda)$ satisfying the three term recurrence relation (\ref{eq:6.1}). Moreover the extremal polynomial $\tilde p$ for which (\ref{eq:6.1b}) becomes an equality, is given by
\begin{equation}\label{eq:6.1c}
\tilde p= \sum_{j=1}^n \frac{w_j^{(1,n)}}{j \sqrt{k_{j-1}^{(1)}}} R_j
\end{equation}
where $(w^{(1,n)})^T=(w_1^{(1,n)},\ldots,w_n^{(1,n)})$ is the eigenvector of $\tilde K_{n,1}$ associated to the eigenvalue $\mu_{1,n}$.
\end{theorem}
\noindent {\bf Proof.}\\
Let $p$ be a polynomial of $\mathcal{P}_n$ written in the basis of the $R_i$'s. We look for the polynomials $p$ such that $y^T(K_{n,0}^{(0)} - \lambda K_{n,0}^{(1)})y=0$. Let $\tilde y$ be the vector of $\mathbb{R}^n$ such that $\tilde y^T=(y_1,\ldots,y_n)$.
\begin{eqnarray*}
y^T(K_{n,0}^{(0)} - \lambda K_{n,0}^{(1)})y&=& y^T \left(\begin{array}{cc}k_0^{(0)}&0\\0& K_{n,1}^{(0)} - \lambda K_{n,1}^{(1)}\end{array} \right) y\\
&=& k_0^{(0)} y_0^2 + \tilde y^T (K_{n,1}^{(0)} - \lambda K_{n,1}^{(1)}) \tilde y\\
&=&k_0^{(0)} y_0^2 + \tilde y^T ((K_{n,1}^{(1)})^{1/2}\left(\tilde K_{n,1} -\lambda I\right)(K_{n,1}^{(1)})^{1/2}) \tilde y
\end{eqnarray*}
Therefore, if $w^{(1,n)}=(K_{n,1}^{(1)})^{1/2}) \tilde y$ is the eigenvector of $\tilde K_{n,1}$ associated to the eigenvalue $\mu_{1,n}$ and if $y_0=0$, we have $y^T(K_{n,0}^{(0)} - \lambda K_{n,0}^{(1)})y=0$. The entries $y_j$, for $j=1,\ldots,n$, are equal to $\frac{w_j^{(1,n)}}{j \sqrt{k_{j-1}^{(1)}}}$. Hence (\ref{eq:6.1c}) holds
$\square$

In some cases an explicit upper bound of this constant can be obtained by using the Newton method on the polynomial $A_n (\lambda)$ (see \cite{DEH2}, \cite{DRA}). Indeed we have $0<x_1=-\frac{A_n (0)}{A_n^\prime (0)}< \mu_{1,n}$. $A_n (0)$ is easily obtained from the Cholesky decomposition of $\tilde K_{n,1}$. $\tilde K_{n,1}= \tilde G_{n,1} ( \tilde G_{n,1})^T$ with
\begin{eqnarray*}
\tilde G_{n,1}(i,i)&=&\frac{1}{i}\sqrt{\frac{k_{i}^{(0)}}{k_{i-1}^{(1)}}}, \quad  i = 1,\ldots,n,\\
\tilde G_{n,1}(i,i-1)&=& - \frac{\sigma_{i-1}}{i-1} \sqrt{\frac{k_{i-1}^{(0)}}{k_{i-1}^{(1)}}},  \quad  i = 2,\ldots,n.
\end{eqnarray*}
Since $A_n(\lambda)=\det(\lambda I - \tilde K_{n,1})$, we have
\begin{equation}\label{eq:6a.1}
A_n(0) = (-1)^n(\det(\tilde G_{n,1}))^2 = (-1)^n\prod_{i=1}^n (\tilde G_{n,1}(i,i))^2= (-1)^n \frac{1}{(n!)^2}\prod_{i=1}^n \frac{k_{i}^{(0)}}{k_{i-1}^{(1)}}.
\end{equation}
The expression of $A_n^\prime (0)$ is obtained from
\begin{equation}\label{eq:6a.2}
A_n^\prime(0)=A_{n-1} (0) -( \frac{k_{n}^{(0)}}{n^2 k_{n-1}^{(1)}}+ \frac{\sigma_{n-1}^2 }{({n-1})^2}\frac{ k_{n-1}^{(0)}}{k_{n-1}^{(1)}} )A_{n-1}^\prime (0)- \frac{\sigma_{n-1}^2 (k_{n-1}^{(0)})^2}{(n-1)^4 k_{n-1}^{(1)}k_{n-2}^{(1)}} A_{n-2}^\prime (0)
\end{equation}
with $A_0^\prime(0)=0$ and $A_1^\prime(0)=1$.\\
Moreover the $qd$ algorithm, applied to the sequence $\{A_j(\lambda)\}_{j=1}^n$, also provides in some cases an explicit lower bound of the Markov-Bernstein constant (see \cite{DRA}). In this case the starting relations of this algorithm are
\begin{eqnarray}
q_{j}^{(0)}&=&\frac{k_{j}^{(0)}}{j^2 k_{j-1}^{(1)}},\quad j=1,\ldots,n, \label{eq:6b.1}\\
e_{j-1}^{(0)}&=&  \frac{\sigma_{j-1}^2 }{({j-1})^2}\frac{ k_{j-1}^{(0)}}{k_{j-1}^{(1)}}, \quad j=2,\ldots,n.\label{eq:6b.2}
\end{eqnarray}
Then, next sequences $\{q_{j}^{(r)}\}$ and $\{e_{j}^{(r)}\}$ are computed as follows, for $r=0,1, \ldots$
\begin{eqnarray}
e_{0}^{(r)}&=&e_{n}^{(r)}=0, \nonumber\\
q_{j+1}^{(r)}+e_{j+1}^{(r)}&=& q_{j+1}^{(r+1)}+e_{j}^{(r+1)}, \quad j=0,\ldots,n-1,\nonumber\\
q_{j+1}^{(r)}e_{j}^{(r)}&=&q_{j}^{(r+1)}e_{j}^{(r+1)}, \quad j=1,\ldots,n-1 \label{eq:6b.3}
\end{eqnarray}
and the following inequalities hold (see \cite{DRA} Corollary 2.4).
\begin{equation}\label{eq:6b.4}
\mu_{1,n} < q_{n}^{(r)}< q_{n}^{(r-1)}<\cdots< q_{n}^{(0)}.
\end{equation}
Another way to obtain some informations about the zeros of $A_n(\lambda)$ is to use Blumenthal Theorem that we recall below in its simplified version (see Chihara \cite{CHI} page 122, Theorem 4.1).
%theorem
\begin{theorem}\label{1.5}
Let $A_n(\lambda)=(\lambda-B_n)A_{n-1}(\lambda)- C_n A_{n-2}(\lambda)$ be a three term recurrence relation in which
$$
\lim_{n \to \infty} B_n= \nu \quad \mbox{ and } \lim_{n \to \infty} C_n = \eta >0.
$$
Let us set $\sigma = \nu - 2 \sqrt{\eta}$ and $\tau = \nu + 2 \sqrt{\eta}$.\\
Then, the set $X$ of all the zeros of all the polynomials $A_n(\lambda)$ is dense in $[\sigma,\tau]$.
\end{theorem}
Some isolated zeros can be lain on the outside of $[\sigma,\tau]$.
%Remark
\begin{remark}\label{1.6}
If $\sigma=0$ and $\eta>0$, or if $\sigma=\eta=0$, since all the zeros of $A_n(\lambda)$ are positive, we can conclude that $\lim_{n \to \infty} \mu_{1,n}=0$.
\end{remark}

Before studying the different cases of coherent pairs, some properties of monic Laguerre-Sonin and Jacobi polynomials are recalled.\\
%Monic Laguerre-Sonin polynomials
{\bf Monic Laguerre-Sonin polynomials $L_n^{\alpha}(x)$, $\alpha >-1$.}\\
They satisfy the following three term recurrence relation (see Chihara \cite{CHI} Relation (2.30) page 154):
\begin{equation}\label{eq:6.2}
L_{n+1}^{\alpha}(x)=(x-2n-\alpha-1)L_n^{\alpha}(x)-n(n+\alpha)L_{n-1}^{\alpha}(x), \quad \forall n \ge 1
\end{equation}
with $L_0^{\alpha}(x)=1$ and $L_1^{\alpha}(x)=x-\alpha-1$.\\
Their norm is 
\begin{equation}\label{eq:6.3}
k_n^{\alpha}=\int_0^\infty (L_n^{\alpha}(x))^2 x^\alpha e^{-x} dx=n! \Gamma(\alpha +n+1)
\end{equation}
 where $\Gamma$ is the Gamma function.\\
The explicit formula  of $L_n^{\alpha}(x)$ is
\begin{equation}\label{eq:6.4}
L_n^\alpha (x)=(-1)^n n! \sum_{m=0}^n \frac{(-1)^m}{m!} \left(\begin{array}{c}n+\alpha\\n-m \end{array} \right)x^m
\end{equation}
where $\left(\begin{array}{c}a\\n \end{array} \right)$ denotes the binomial coefficient.\\
Moreover 
\begin{equation}\label{eq:6.5}
\frac{d}{dx} L_n^\alpha (x) = n L_{n-1}^{\alpha+1} (x).
\end{equation}
%Monic Jacobi polynomials
{\bf Monic Jacobi polynomials $P_n^{\alpha,\beta}(x)$, $\alpha >-1$, $\beta >-1$.}\\
They satisfy the following three term recurrence relation (see Chihara \cite{CHI} Relation (2.29) page 153):
\begin{eqnarray}
P_{n+1}^{(\alpha,\beta)}(x)&=& (x-\frac{\beta^2-\alpha^2}{(2n+\alpha +\beta)(2n+\alpha +\beta+2)})P_n^{(\alpha,\beta)}(x) \nonumber\\
&&- \frac{4 n (n+\alpha)(n+\beta)(n+\alpha +\beta)}{(2n+\alpha +\beta)^2(2n+\alpha +\beta+1)(2n+\alpha +\beta-1)} P_{n-1}^{(\alpha,\beta)}(x) \label{eq:7.1}
\end{eqnarray}
with $P_{0}^{(\alpha,\beta)}(x)=1$ and $P_{1}^{(\alpha,\beta)}(x)=x-\frac{\beta-\alpha}{\alpha +\beta +2}$.\\
Their norm is 
\begin{eqnarray}\label{eq:7.2}
k_n^{(\alpha,\beta)}&=&\int_{-1}^1 (P_n^{(\alpha,\beta)}(x))^2 (1-x)^{\alpha}(1+x)^{\beta} dx \nonumber\\
&=&2^{2n+\alpha + \beta+1} n! \frac{\Gamma(n+\alpha +1)\Gamma(n+\beta +1)\Gamma(n+\alpha +\beta +1)}{(2n+\alpha +\beta +1)(\Gamma(2n+\alpha +\beta +1))^2}.
\end{eqnarray}
Note that (\ref{eq:7.2}) could also be written for more convenience as
\begin{equation}\label{eq:7.6}
P_{n+1}^{(\alpha,\beta)}(x)= (x-\frac{\beta^2-\alpha^2}{(2n+\alpha +\beta)(2n+\alpha +\beta+2)})P_n^{(\alpha,\beta)}(x) - \frac{k_n^{(\alpha,\beta)}}{k_{n-1}^{(\alpha,\beta)}} P_{n-1}^{(\alpha,\beta)}(x).
\end{equation}
We also have two other recurrence relations
\begin{eqnarray}
 P_n^{(\alpha,\beta)} (x)  
&=&   P_n^{(\alpha+1,\beta)} (x) - \frac{2 n (n+\beta)}
{ (2 n +\alpha+\beta)(2 n +\alpha+\beta+1)}  P_{n-1}^{(\alpha+1,\beta)} (x),\label{eq:L8}\\
P_n^{(\alpha,\beta)} (x) 
&=&   P_n^{(\alpha,\beta+1)} (x)+  \frac{2 n (n+\alpha)}
{ (2 n +\alpha+\beta)(2 n +\alpha+\beta+1)}  P_{n-1}^{(\alpha,\beta+1)} (x),\label{eq:L9}
\end{eqnarray}
Two explicit formulas will be used
\begin{eqnarray}
P_n^{(\alpha,\beta)}(x)&=& 2^n \frac{\Gamma(\alpha +n+1)}{\Gamma(\alpha+\beta +2n+1)}\sum_{m=0}^n \left(\begin{array}{c}n\\m\end{array}\right)\frac{\Gamma(\alpha+\beta +n+m+1)}{\Gamma(\alpha +m+1)}(\frac{x-1}{2})^m \label{eq:7.3}\\
&=& (-2)^n \frac{\Gamma(\beta +n+1)}{\Gamma(\alpha+\beta +2n+1)}\sum_{m=0}^n  \left(\begin{array}{c}n\\m\end{array}\right)\frac{\Gamma(\alpha+\beta +n+m+1)}{\Gamma(\beta +m+1)}(-\frac{x+1}{2})^m.\label{eq:7.4}
\end{eqnarray}
Moreover
\begin{equation}\label{eq:7.5}
\frac{d}{dx}P_n^{(\alpha,\beta)}(x) = n P_{n-1}^{(\alpha+1,\beta+1)}(x).
\end{equation}

%**************************************************************************************************************************

%_______________________________________________section 2 ..........................

%**************************************************************************************************************************

\section{Laguerre case} 

%**************************************************************************************************************************

%____________________________________________section 2.1________________________________

%**************************************************************************************************************************
\subsection{Laguerre case (a)} 

The measure associated to $c_0$ is $(x-\xi)x^{\alpha -1} e^{-x}dx$ with $\alpha >0$ and $\xi < 0$, and the one associated to
 $c_1$ is the Laguerre-Sonin measure $x^{\alpha } e^{-x}dx$. Therefore $T_n(x)= L_n^\alpha(x)$. The monic polynomials, orthogonal with respect to $c_1$, will be denoted by $ P_n(x,\xi)$. Their expression, depending on the monic Laguerre polynomials $L_j^{\alpha-1}(x)$, has already been given in the thesis of P\'erez \cite{P} (Corollary 5.3.3 and Theorem 5.3.4). Nevertheless we give a direct proof of this expression. In this case we also obtain the square norm $k_n^{(0)}$ and $\sigma_n$.
%theorem
\begin{theorem}\label{2.1}
The monic polynomials $P_n(x,\xi)$, orthogonal with respect to $c_1$, satisfy the following relation:
\begin{equation}\label{eq:8.1}
P_n(x,\xi)= N_n^{\alpha-1}(x,\xi) \frac{k_n^{\alpha-1}}{L_n^{\alpha-1}(\xi)}
\end{equation}
where $N_n^{\alpha-1}(x,\xi)$ is the reproducing kernel associated to the measure $x^{\alpha -1} e^{-x}dx$.
\begin{equation*}
N_n^{\alpha-1}(x,\xi) =\sum_{j=0}^n \frac{L_j^{\alpha-1}(x)L_j^{\alpha-1}(\xi)}{k_j^{\alpha-1}}. 
%\label{eq:8.2}
\end{equation*}
The square norm $k_n^{(0)}$ of $P_n(x,\xi)$ and $\sigma_n$ are given by
\begin{eqnarray}
k_n^{(0)}&=&-n! \Gamma(\alpha +n) \frac{L_{n+1}^{\alpha-1}(\xi)}{L_{n}^{\alpha-1}(\xi)}, \quad \forall n \ge 0, \label{eq:8.6}\\
\sigma_n&=&n(n+\alpha)\frac{L_{n}^{\alpha-1}(\xi)}{L_{n+1}^{\alpha-1}(\xi)}, \quad \forall n \ge 1.\label{eq:9.3}
\end{eqnarray}
\end{theorem}
\noindent {\bf Proof.}\\
Let us write $(x-\xi) P_n(x,\xi)$ in the basis of monic Laguerre polynomials $L_j^{\alpha-1}(x)$.
\begin{equation}\label{eq:8.3}
(x-\xi) P_n(x,\xi)= \sum_{j=0}^{n+1} \theta_{j,n+1} L_j^{\alpha-1}(x), \quad \theta_{n+1,n+1}=1.
\end{equation}
Let us denote by $c^{\alpha-1}$ the linear functional defined from the measure $x^{\alpha-1} e^{-x}$ on the support $\Omega=]0,+\infty[$. We have.
\begin{eqnarray*}
c^{\alpha-1}(L_i^{\alpha-1}(x)(x-\xi) P_n(x,\xi))&=& c_0 (L_i^{\alpha-1}(x) P_n(x,\xi)) =0 \mbox{ for } i=0,\ldots,n-1,\\
&=& \sum_{j=0}^{i} \theta_{j,n+1} c^{\alpha-1}(L_i^{\alpha-1}(x) L_j^{\alpha-1}(x))\\
&=&  \theta_{i,n+1} k_i^{\alpha-1}.
\end{eqnarray*}
Therefore $\theta_{i,n+1}=0$ for $i=0\ldots,n-1$, and 
\begin{equation}\label{eq:8.4}
\theta_{n,n+1}=\frac{k_n^{(0)}}{k_n^{\alpha-1}}.
\end{equation}
We put $x=\xi$ in (\ref{eq:8.3}). Then we obtain
\begin{equation}\label{eq:8.5}
\theta_{n,n+1}=-\frac{L_{n+1}^{\alpha-1}(\xi)}{L_{n}^{\alpha-1}(\xi)}.
\end{equation}
The expression (\ref{eq:8.6}) of the square norm $k_n^{(0)}$ is obtained from (\ref{eq:8.4}) and (\ref{eq:8.5}).\\
From (\ref{eq:8.3}) we have
\begin{equation}\label{eq:9.1}
P_n(x,\xi)=\left(L_{n+1}^{\alpha-1}(x)- \frac{L_{n+1}^{\alpha-1}(\xi)}{L_{n}^{\alpha-1}(\xi)}L_{n}^{\alpha-1}(x)\right)\frac{1}{x-\xi}.
\end{equation}
(\ref{eq:8.1}) is obtained by using the definition of the reproducing kernel $N_n^{\alpha-1}(x,\xi)$ associated to the Christoffel-Darboux formula.\\
Now let us determine explicitly the sequence $\{\sigma_n\}_{n\ge 1}$ as a function depending on $\xi$. From (\ref{eq:1.1}), (\ref{eq:6.5}) and (\ref{eq:8.1}) we have
\begin{eqnarray}
\sigma_n&=&\frac{n}{n+1} \frac{P_{n+1}^\prime(x,\xi)-(n+1) L_{n}^{\alpha}(x)}{P_{n}^\prime(x,\xi)} \nonumber\\
&=& \frac{n}{n+1}\frac{1}{P_{n}^\prime(x,\xi)} \sum_{j=1}^{n} j L_{j-1}^{\alpha}(x) \frac{L_j^{\alpha-1}(\xi)}{L_{n+1}^{\alpha-1}(\xi) }\frac{k_{n+1}^{\alpha-1}}{k_j^{\alpha-1}}  \label{eq:9.2}
\end{eqnarray}
To obtain $\sigma_n$ it is sufficient to give the ratio of the coefficients of $x^{n-1}$ in the numerator and the denominator of (\ref{eq:9.2}).  Hence (\ref{eq:9.3}) holds.$\square$\\
Note that (\ref{eq:9.3}) is also given in \cite{MP}.\\
At last the three term recurrence relation satisfied by the polynomials $A_n(\lambda)$ is given.
%theorem
\begin{theorem}\label{2.2}
The following Markov-Bernstein inequality holds:
$$
c_1((p^\prime)^2) \le \frac{1}{\mu_{1,n}} c_0(p^2), \quad \forall p \in \mathcal{P}_n
$$
where $\mu_{1,n}$ is the smallest zero of the polynomials $A_n(\lambda)$ satisfying the following three term recurrence relation
\begin{equation}\label{eq:10.1}
A_n(\lambda)=(\lambda -2+\frac{\xi-\alpha}{n})A_{n-1}(\lambda)-(1+\frac{\alpha}{n-1})A_{n-2}(\lambda), \quad \forall n\ge 2.
\end{equation}
$A_0(\lambda)=1$ and $A_1(\lambda)= \lambda -\alpha-1+\xi - \frac{\xi}{\xi-\alpha}$.
\end{theorem}
\noindent {\bf Proof.}\\
We use the expressions of $k_n^{(0)}$ and $\sigma_n$ in (\ref{eq:6.1}). We obtain
$$
A_n(\lambda)=(\lambda +\frac{1}{n}\frac{L_{n+1}^{\alpha-1}(\xi)}{L_{n}^{\alpha-1}(\xi)}+(n+\alpha-1)\frac{L_{n-1}^{\alpha-1}(\xi)}{L_{n}^{\alpha-1}(\xi)})A_{n-1}(\lambda)-\frac{n+\alpha-1}{n-1}A_{n-2}(\lambda), \quad \forall n\ge 2.
$$
This relation is simplified by using the three term recurrence relation (\ref{eq:6.2}). Hence (\ref{eq:10.1}) holds.
$\square$
%Remark
\begin{remark}\label{2.3}
In the case of Markov-Bernstein inequalities associated to a bilinear functional $\tilde a_\lambda(p,q)=c_1(pq)+ \lambda c_1(p^\prime q^\prime)$ for $p,q \in \mathcal{P}$, the Markov-Bernstein constant $1/\sqrt{\mu_{1,n}}$ is obtained with the smallest zero $\mu_{1,n}$ of the polynomial $\tilde A_n(\lambda)$ satisfying the following three term recurrence relation (see Milovanovi\'c et al \cite{MMR} page 576, Draux and Elhami \cite{DEH2} Theorem 6.3)
\begin{equation}\label{eq:11.1}
\tilde A_n(\lambda)=(\lambda -2-\frac{\alpha}{n})\tilde A_{n-1}(\lambda)-(1+\frac{\alpha}{n-1})\tilde A_{n-2}(\lambda), \quad \forall n\ge 2,
\end{equation}
$\tilde A_0(\lambda)=1$ and $\tilde A_1(\lambda)=\lambda -\alpha -1$.\\
(\ref{eq:10.1}) tends to (\ref{eq:11.1}) when $\xi$ tends to $0$. $\xi$ can be considered as a perturbation of the classical case $\tilde a_\lambda(p,q)$.\\
Moreover the following three term recurrence relation
$$
p_{n}(\lambda)=(\lambda-2)p_{n-1}(\lambda)-p_{n-2}(\lambda)
$$
can be transformed as
$$
p_{n}^*(y)=2yp_{n-1}^*(y)-p_{n-2}^*(y)
$$
by using $p_{n}^*(y)=p_n(\frac{\lambda-2}{2})$. We obtain the three term recurrence relation corresponding to Chebyshev polynomials. Therefore Relations (\ref{eq:10.1}) and (\ref{eq:11.1}) can also be considered as a perturbed Chebyshev  three term recurrence relation by a $\ell^2$ perturbation of the coefficients. Such perturbations of coefficients of a three term recurrence relation appear, for example, in relations studied in \cite{ADT}.
\end{remark}
%Remark
\begin{remark}\label{2.3aa}
The fact that the polynomials $A_n(\lambda)$ satisfying Relation (\ref{eq:11.1}) are monic co-recursive Pollaczek polynomials, was already established by P\'erez in her doctoral dissertation (see \cite{P} page 28; see also Marcell\'an et al \cite{MPP2}). Monic co-recursive Pollaczek polynomials satisfy the three term recurrence relation of monic Pollaczek polynomials $P^\gamma_n(x;a,b,c)$ (see Chihara \cite{CHI} page 185 Relation 5.9), but the polynomial of degree 1 is different.
\begin{eqnarray}\label{eq:n.1}
P^\gamma_n(x;a,b,c)
&= &(x+\frac{b}{n+\gamma+a+c-1})P^\gamma_{n-1}(x;a,b,c) \nonumber\\
&&-\frac{(n+c-1)(n+2\gamma+c-2)}{(n+\gamma+a+c-1)(n+\gamma+a+c-2)}P^\gamma_{n-2}(x;a,b,c),\\
P^\gamma_1(x;a,b,c)
&= &x+\frac{b}{\gamma+a+c}. \nonumber
\end{eqnarray}
Indeed, if $A^*_n(y)=A_n(\frac{\lambda-2}{2})$, we have 
\begin{eqnarray}\label{eq:n.2}
A^*_n(y)&=& (y+\frac{\xi-\alpha}{2n}) A^*_{n-1}(y)-\frac{1}{4}(1+\frac{\alpha}{n-1})  A^*_{n-2}(y),\\
A^*_1(y)&=&y+ \frac{\xi-\alpha}{2}+ \frac{\alpha}{2(\xi-\alpha)}. \nonumber
\end{eqnarray}
Relations (\ref{eq:n.1}) and (\ref{eq:n.2}) are identical if and only if
$$
\gamma+a+c-1=0, \quad b= \frac{\xi-\alpha}{2}, \quad \alpha =2(\gamma+c-1) \mbox{ and } (c-1)(2\gamma+c-2)=0.
$$
Therefore, either $c=1$ or $c=2(1-\gamma)$.\\
If $c=1$, then $\alpha =2 \gamma$, $a=-\gamma$ and $b= \frac{\xi-\alpha}{2}$.\\
If $c=2(1-\gamma)$, then $\alpha =c=2 (1-\gamma)$, $a=\gamma-1$ and $b= \frac{\xi-\alpha}{2}$.\\
Hence the corresponding Pollaczek polynomials are $P^\gamma_n(x;-\gamma,b,1)$ or $P^\gamma_n(x;\gamma-1,b,2 (1-\gamma))$ with $b= \frac{\xi-\alpha}{2}$ (see also \cite{ADT} for a use of these two families of polynomials). $A^*_1(y)$ is different from $P^\gamma_1(x;a,b,c) = x+\frac{b}{\gamma+a+c}=x+b= x+ \frac{\xi-\alpha}{2}$.
\end{remark}
%property
\begin{property}\label{2.3a}
$$
\lim_{n \to \infty} \mu_{1,n}=0.
$$
\end{property}
\noindent {\bf Proof.}\\
We use Blumenthal Theorem \ref{1.5}. $\lim_{n \to \infty} B_n = 2$ and $\lim_{n \to \infty} C_n =1$. Therefore $\sigma=0$ and $\tau =4$. Hence the property holds by using Remark \ref{1.6}.
$\square$

%**************************************************************************************************************************

%____________________________________________section 2.2____________________________

%**************************************************************************************************************************
\subsection{Laguerre case (b)} 

The measure associated to $c_0$ is $e^{-x}dx + M \delta(0)$ with $M\ge0$, and the one associated to
 $c_1$ is the Laguerre measure $e^{-x}dx$. Therefore $T_n(x)$ is the monic classical Laguerre polynomial $L_n(x)$. Let us denote by $P_n(x,M)$, $n\ge0$, the monic polynomials orthogonal with respect to $c_0$.\\
 In order to obtain an explicit expression of $P_n(x,M)$, it will be written in the basis of monic Laguerre polynomial $L_j(x)$.
\begin{equation}\label{eq:11.2}
P_n(x,M)= \sum_{j=0}^n \theta_{j,n} L_j(x), \quad \theta_{n,n}=1.
\end{equation}
%property
\begin{property}\label{2.6}
The explicit expressions of $P(x,M)$ in the basis of monic Laguerre polynomials $L_j(x)$, of its square norm $k_n^{(0)}$ and of $\sigma_n$ are
\begin{eqnarray}
P_n(x,M)&=& L_n(x)+\sum_{j=0}^{n-1} (-1)^{n+1+j} \frac{n! M}{j! (nM+1)} L_j(x), \label{eq:14.1} \\
k_n^{(0)}&=& (n!)^2  \frac{ (n+1)M+1}{nM+1}, \quad \forall n \ge 0, \label{eq:14.2}\\
\sigma_n &=& - \frac{ n(1+nM)}{(n+1)M+1}, \quad \forall n \ge 1. \label{eq:14.3}
\end{eqnarray}
\end{property}
\noindent {\bf Proof.}\\
To obtain the $\theta_{j,n}$'s we use the orthogonality relations of $P_n(x,M)$ for $i=0,\ldots,n-1$.
\begin{eqnarray*}
c_0(L_i(x) P_n(x,M)) &=& \sum_{j=0}^{n} \theta_{j,n} \left( c_1(L_i(x) L_j(x)) + M (L_i(0) L_j(0))\right) =0\\
&=& \theta_{i,n} k_i^{(1)} + M (-1)^i i!  \sum_{j=0}^{n} \theta_{j,n} (-1)^j j! 
\end{eqnarray*}
Thus 
\begin{equation}\label{eq:l.1}
\theta_{i,n}= \frac{(-1)^{i+1}}{ i! } M \sum_{j=0}^{n} \theta_{j,n} (-1)^j j! \quad \mbox{ for }i=0,\ldots,n-1.
\end{equation}
Let us set $K=\sum_{j=0}^{n} \theta_{j,n} (-1)^j j! $. Therefore
$$
K= (-1)^n n!+ \sum_{j=0}^{n-1} \theta_{j,n} (-1)^j j! = (-1)^n n! +  \sum_{j=0}^{n-1} (-1)^{2j+1} M K 
$$
by using (\ref{eq:l.1}) in the definition of $K$. Hence $K=(-1)^{n} \frac{n! }{ nM+1}$ and
\begin{eqnarray*}
\theta_{j,n} &=& (-1)^{n+1+j} \frac{n! M}{j! (nM+1)}, \mbox{ for } j=0,\ldots,n-1,\\
\theta_{n,n} &=&1.
\end{eqnarray*}
Hence (\ref{eq:14.1}) holds.\\
Now we prove (\ref{eq:14.2}).
\begin{eqnarray*}
c_0((P_n(x,M))^2)&=&c_1((P_n(x,M))^2) + M (P_n(0,M))^2\\
&=& k_n^{(1)} + \sum_{j=0}^{n-1} ( \frac{n! M}{j! (nM+1)})^2 k_j^{(1)} + M\left(L_n(0)+\sum_{j=0}^{n-1} (-1)^{n+1+j} \frac{n! M  L_j(0)}{j! (nM+1)}\right)^2\\
&=&(n!)^2  \frac{ (n+1)M+1}{nM+1}.
\end{eqnarray*}
At last $\sigma_n$ is obtained from (\ref{eq:1.1}).
$$
\sigma_n=\frac{n}{n+1} \frac{P_{n+1}^\prime(x,M)-(n+1) L_n(x)}{P_{n}^\prime(x,M)}.
$$
It is sufficient to give the ratio of the coefficients of $x^{n+1}$ of the numerator and of the denominator.
$$
P_{n+1}^\prime(x,M)=(n+1) L_n^1(x)+\sum_{j=1}^{n} (-1)^{n+2+j} \frac{(n+1)! M}{j! ((n+1)M+1)} j L_{j-1}^1(x).
$$
Its coefficient of $x^{n+1}$ is $-n(n+1)\frac{ n(n+2)M+n+1}{(n+1)M+1}$. The one of $(n+1) L_n(x)$ is $-n^2(n+1)$. Hence (\ref{eq:14.3}) holds.
$\square$\\
If we rather use the following orthogonality relations $c_0(x^i P_n(x,M))=0$ for $i=0,\ldots,n-1$, two kinds of relations are obtained.\\
For $i=0$
\begin{equation*}
%\label{eq:14.4}
\theta_{0,n}+ M  \sum_{j=0}^n \theta_{j,n} L_j(0)= \theta_{0,n}+ M  \sum_{j=0}^n (-1)^j \theta_{j,n} j!=0.
\end{equation*}
For $i=1,\ldots,n-1$, 
$$
\sum_{j=0}^i \theta_{j,n} c_1(x^i L_j(x))=0.
$$
These last relations need to compute $c_1(x^i L_j(x))$ for $i\ge j$. It is much more complicated, but the expression of $c_1(x^i L_j(x))$ for $i\ge j$ is new. Moreover an original identity with a sum of products of binomial coefficients is deduced from this computation. It is the reason for which we give these results in the following Lemma and Corollary.
%lemma
\begin{lemma}\label{2.4}
\begin{equation}\label{eq:12.1}
\forall i \ge j, \quad c_1(x^i L_j(x))= i! j! \left(\begin{array}{c}i\\j \end{array} \right).
\end{equation}
\end{lemma}
\noindent {\bf Proof.}\\
By using the formal expression of $L_n(x)$ given by (\ref{eq:6.4}) we have
$$
c_1(x^i L_j(x))=(-1)^j j! \sum_{m=0}^j \frac{(-1)^m}{m!}\left(\begin{array}{c}j\\m \end{array} \right)c_1(x^{m+i}).
$$
From the definition of the Gamma function we have
\begin{equation}\label{eq:12.2}
c_1(x^{m+ i})= \Gamma(m+i+1)=(m+i)!.
\end{equation}
Therefore
\begin{equation}\label{eq:12.3}
c_1(x^i L_j(x))= (-1)^j j! i! \sum_{m=0}^j (-1)^m \left(\begin{array}{c}j\\m \end{array} \right) \left(\begin{array}{c}m+i\\i \end{array} \right).
\end{equation}
For $i=j$, (\ref{eq:12.1}) corresponds to the square norm $k_j^{(1)}=(j!)^2$.\\
To prove (\ref{eq:12.1}) for $i> j$, we use the three term recurrence relation satisfied by the polynomials $L_j(x)$
\begin{equation}\label{eq:12.4}
L_{j+1}(x)=(x-2j-1) L_j(x) - j^2 L_{j-1}(x), \quad \forall j\ge 0.
\end{equation}
This relation is multiplied by $x^j$ and $c_1$ is applied to the result. We obtain
\begin{equation*}
%\label{eq:12.5}
c_1(x^{j+1} L_j(x))=(2j+1)! (j!)^2 + j^2 c_1(x^{j}L_{j-1}(x)), \quad \forall j\ge 0.
\end{equation*}
For $j=1$ we have $c_1(x^{j+1} L_j(x))=4$. Therefore (\ref{eq:12.1}) holds for $i=j+1=2$. Then, by induction
\begin{eqnarray*}
c_1(x^{j+1} L_j(x))&=&(2j+1)(j!)^2 + j^2 j! (j-1)!  \left(\begin{array}{c}j\\j-1 \end{array} \right)\\
&=&j! (j+1)!  \left(\begin{array}{c}j+1\\j \end{array} \right).
\end{eqnarray*}
Therefore (\ref{eq:12.1}) holds for $i=j+1$, $\forall j \ge 0$.\\
Now (\ref{eq:12.1}) is assumed to hold for $\ell \in \mathbb{N}$, $\ell$ fixed such that $i=j+\ell$, $\forall j \ge 0$.
$$
c_1(x^{j+\ell} L_j(x))=j! (j+\ell)!  \left(\begin{array}{c}j+\ell\\j \end{array} \right), \quad \forall j\ge 0.
$$
From (\ref{eq:12.4}) we have
\begin{eqnarray*}
c_1(x^{j+\ell+1} L_j(x))&=&c_1(x^{j+\ell} L_{j+1}(x)) +(2j+1) c_1(x^{j+\ell} L_j(x)) + j^2 c_1(x^{j+\ell} L_{j-1}(x))\\
&=&(2j+\ell+1)j! (j+\ell)! \left(\begin{array}{c}j+\ell\\j \end{array} \right) + j^2 c_1(x^{j+\ell} L_{j-1}(x)).
\end{eqnarray*}
$c_1(x^{\ell+1} L_{0}(x))=((\ell+1)!)^2$, therefore (\ref{eq:12.1}) holds for $j=0$ and $i=\ell+1$. After that, by induction
\begin{eqnarray*}
c_1(x^{j+\ell+1} L_j(x)) &=&(2j+\ell+1)j! (j+\ell)! \left(\begin{array}{c}j+\ell\\j \end{array} \right) + j^2 (j-1)! (j+\ell)! \left(\begin{array}{c}j+\ell\\j-1 \end{array} \right) \\
&=& \frac{((j+\ell+1)!)^2}{(\ell+1)!} = j! (j+\ell+1)! \left(\begin{array}{c}j+\ell+1\\j \end{array} \right) .
\end{eqnarray*}
Hence (\ref{eq:12.1}) holds.
$\square$\\
From (\ref{eq:12.3}) the following corollary is deduced:
%corollary
\begin{corollary}\label{2.5}
$$
 (-1)^j  \sum_{m=0}^j (-1)^m \left(\begin{array}{c}j\\m \end{array} \right) \left(\begin{array}{c}m+i\\i \end{array} \right)=\left(\begin{array}{c}i\\j \end{array} \right), \quad \forall i \ge j.
 $$
\end{corollary}
Finally, by using (\ref{eq:14.2}) and (\ref{eq:14.3}) in (\ref{eq:6.1}) we have the following theorem concerning the Markov-Bernstein inequalities.
%theorem
\begin{theorem}\label{2.7}
The following Markov-Bernstein inequality holds:
$$
c_1((p^\prime)^2) \le \frac{1}{\mu_{1,n}} c_0(p^2), \quad \forall p \in \mathcal{P}_n
$$
where $\mu_{1,n}$ is the smallest zero of the polynomials $A_n(\lambda)$ satisfying the following three term recurrence relation
\begin{equation}\label{eq:16.1}
A_n (\lambda) = (\lambda -2)A_{n-1} (\lambda)-A_{n-2} (\lambda), \quad \forall n\ge 2.
\end{equation}
$A_0(\lambda)=1$ and $A_1(\lambda)= \lambda-\frac{ 1+2M}{1+M}$.
\end{theorem}
A lower bound of $\mu_{1,n}$ is obtained by using the Newton method.
%property
\begin{property}\label{2.8}
The smallest zero $\mu_{1,n}$ of $A_n (\lambda)$ is such that
\begin{equation}\label{eq:16.2}
\mu_{1,n}>\frac{1+(n+1)M}{n(n+1)(3+(n+2)M)}, \quad \forall n \ge 1.
\end{equation}
\end{property}
\noindent {\bf Proof.}\\
$A_n(0)$ is obtained from (\ref{eq:6a.1}) and (\ref{eq:14.2}).
$$
A_n(0)=(-1)^n \frac{1+(n+1)M}{1+M}.
$$
$A_n^\prime(0)= (-1)^{n+1}\frac{n(n+1)}{6(1+M)}(3+(n+2)M)$ is obtained from (\ref{eq:6a.2}) by using a proof by recurrence. Hence the result holds.
$\square$\\
Then, the following inequality is deduced from the previous result
%Corollary
\begin{corollary}\label{2.9}
$$
c_1((p^\prime)^2) < \frac{n(n+1)(3+(n+2)M)}{1+(n+1)M} c_0(p^2), \quad \forall p \in \mathcal{P}_n.
$$
\end{corollary}

A better lower bound of $\mu_{1,n}$ is provided by using the Laguerre method (see Draux \cite{DRA}). Unfortunately this formal expression is very complicated. Nevertheless we give it.
$$
0<\tilde x_1=-\frac{n A_n(0)}{A_n^\prime(0)-\sqrt{H_n(0)}}<\mu_{1,n}
$$
with $H_n(0)=(n-1)^2 (A_n^\prime(0))^2-n(n-1)A_n(0)A_n^{\prime\prime}(0)$.\\
$A_n^{\prime\prime}(0))= 2 A_n^\prime(0) -2 A_{n-1}^{\prime\prime}(0))-A_{n-2}^{\prime\prime}(0))$ and from this relation we obtain by recurrence
$$
A_n^{\prime\prime}(0))= (-1)^n (n-1)_4 \frac{5+(n+3)M}{60(1+M)}
$$
where $(a)_n$ is the Pochhammer symbol: $(a)_n=a (a+1)\ldots(a+n-1)$. By definition $(a)_0=1$. Then
\begin{eqnarray*}
\tilde x_1&=& \frac{1+(n+1)M}{\frac{(n+1)(3+(n+2)M)}{6}-(-1)^{n+1}(n-1)\sqrt{K(n,M)}},\\
K(n,M) &=&\frac{(n+1)^2(3+(n+2)M)^2}{36}-\frac{1+(n+1)M}{60}(n+1)(n+2)(5+(n+3)M).
\end{eqnarray*}

At last we give the entries of the $qd$ algorithm in order to obtain an upper bound of $\mu_{1,n}$. Relations (\ref{eq:6b.1}) and (\ref{eq:6b.2}) give the starting entries.
\begin{eqnarray*}
q_{j}^{(0)}&=&\frac{1+(j+1)M}{1+jM}, \quad \mbox{ for }j=1,\ldots,n,\\
e_{j}^{(0)}&=&\frac{1+jM}{1+(j+1)M}, \quad \mbox{ for }j=1,\ldots,n-1.
\end{eqnarray*}

$$
\left\{
\begin{array}{lll}
e_{0}^{(0)}=e_{n}^{(0)}&=&0,\\
q_{j}^{(1)}&=&\displaystyle{\frac{ j+1}{j} \frac{1+jM}{1+(j+1)M} \frac{p_j^{(1)}}{p_{j-1}^{(1)}}}, \quad j=1,\ldots,n-1,\\
e_{j}^{(1)}&=&\displaystyle{\frac{j}{ j+1} \frac{1+(j+2)M}{1+(j+1)M}  \frac{p_{j-1}^{(1)}}{p_j^{(1)}}}, \quad j=1,\ldots,n-1,\\
q_{n}^{(1)}&=&\displaystyle{\frac{(1+nM)(1+(n+1)M)}{n p_{n-1}^{(1)}}}.
\end{array}
\right.
$$
where $p_j^{(1)}= 1+M(2+j)+\frac{M^2}{6}(2+j)(3+2j)$.
$$
\left\{
\begin{array}{lll}
e_{0}^{(1)}=e_{n}^{(1)}&=&0,\\
q_{j}^{(2)}&=&\displaystyle{\frac{( j+2)(2j+3)}{( j+1)(2 j+1)} \frac{p_{j-1}^{(1)}}{p_j^{(1)}} \frac{p_j^{(2)}}{p_{j-1}^{(2)}}}, \quad j=1,\ldots,n-1,\\
e_{j}^{(2)}&=&\displaystyle{\frac{j(2 j+1)}{ (j+1)(2j+3)}   \frac{p_{j+1}^{(1)}}{p_j^{(1)}} \frac{p_{j-1}^{(2)}}{p_j^{(2)}}}, \quad j=1,\ldots,n-2,\\
e_{n-1}^{(2)}&=& \displaystyle{\frac{(n-1)(2 n-1)}{n (n+1)(2n+1)}   \frac{(1+(n+1)M)^2}{p_{n-1}^{(1)}} \frac{p_{n-2}^{(2)}}{p_{n-1}^{(2)}}},\\
q_{n}^{(2)}&=&\displaystyle{\frac{30(1+(n+1)M)}{(n+1)(2n+1) }} \frac{p_{n-1}^{(1)}}{p_{n-1}^{(2)}}
\end{array}
\right.
$$
where $p_j^{(2)}= 5+5M(3+j)+M^2(3+j)(5+2j)+\frac{M^3}{6}(2+j)(3+j)(5+2j)$.\\
It is easy to see that all the entries $q_{j}^{(1)}$, $e_{j}^{(1)}$, $q_{j}^{(2)}$ and $e_{j}^{(2)}$ satisfy the set of relations (\ref{eq:6b.3}).\\
As a consequence of (\ref{eq:6b.4}) we have
%property
\begin{property}\label{2.8a}
The smallest zero $\mu_{1,n}$ of $A_n (\lambda)$ is such that
\begin{equation}\label{eq:16c.2}
\mu_{1,n}<q_{n}^{(2)}, \quad \forall n \ge 1.
\end{equation}
\end{property}
%Corollary
\begin{corollary}\label{2.9d}
The Markov-Bernstein constant $M_n$ satisfies the following inequality
\begin{equation}\label{eq:16c.4}
\frac{1}{\sqrt{q_{n}^{(2)}}}\le M_n\le \sqrt{\frac{n(n+1)(3+(n+2)M)}{1+(n+1)M} }.
\end{equation}
\end{corollary}
%Corollary
\begin{corollary}\label{2.9b}
$$
\mu_{1,n}=O(\frac{1}{n^2})
$$
and the Markov-Bernstein constant has a behavior as a $O(n)$.
\end{corollary}

Below we give a table of some numerical results for the zero $\mu_{1,n}$ and the lower bounds $x_1$ and $\tilde x_1$ provided by the Newton method and the Laguerre method, as well as the upper bound $q_{n}^{(2)}$ provided by the $qd$ algorithm, for different values of $M$.

\begin{center}
\begin{tabular}{|c|c|c|c|c|c|}
\hline
$M$  & $n$ & $x_1$ & $\tilde x_1$ &$\mu_{1,n}$ & $q_{n}^{(2)}$\\
\hline
 &  20&0.002095238 &0.019766736 &
0.020393972
&0.029017408\\
$1$  &  50&0.000370766 &0.003519638&
0.003649401
&0.005391291 \\
    &  100&0.000096181 & 0.000913322 &
0.000948578
&0.001420806 \\
    &  500&0.000003968 & 0.000037658 &
0.000039164
&0.000059345 \\
\hline
 &   20&0.002233459& 0.021242255 &
0.021922153
&0.031139285 \\
$5$  &  50&0.000381719 &0.003629865 &
0.003763805
&0.005558176 \\
    &  100&0.000097658 & 0.000927797&
0.000963616
&0.001443177 \\
    &  500&0.000003980 & 0.000037778 &
0.000039289
&0.000059534 \\
\hline
 &  20&0.002252829 & 0.021442114 &
0.022128520
&0.031423693 \\
$10$  &  50&0.000383158 &0.003644061 &
0.003778527
&0.005579620 \\
    &  100&0.000097848 &0.000929632 &
0.000965523
&0.001446012 \\
    &  500&0.000003982 & 0.000037793 &
0.000039305
&0.000059558 \\
\hline
 &   20&0.002268704 & 0.021604487 &
0.022296105
&0.031654366 \\
$50$  &  50&0.000384322 &0.003655484 &
0.003790372
&0.005596869 \\
    &  100&0.000098000 &0.000931105 &
0.000967052
&0.001448286 \\
    &  500&0.000003983 & 0.000037805 &
0.000039317
&0.000059577 \\
\hline

\end{tabular}
\end{center}

%**************************************************************************************************************************

%____________________________________________section 2.3_____________________________

%**************************************************************************************************************************
\subsection{Laguerre case (c)} 

The measure associated to $c_0$ is the Laguerre-Sonin measure $x^\alpha e^{-x}dx$ with $\alpha >-1$, and the one associated to
 $c_1$ is  $\dfrac{x^{\alpha+1}}{x-\xi} e^{-x}dx+M \delta(\xi)$ with $\xi \le 0$ and $M \ge 0$. Therefore the monic polynomials $P_n(x)$,  orthogonal with respect to $c_0$, are  the monic Laguerre-Sonin polynomials $L_n^{\alpha}(x)$. Since we will also use the Laguerre-Sonin polynomials $L_n^{\alpha+1}(x)$, we prefer to denote by $c^{\alpha}$ the linear functional $c_0$. Let us denote by $T_n(x,\xi,M)$, $n\ge0$, the monic polynomials orthogonal with respect to $c_1$. Moreover we will denote by $\tilde c_1$ the linear functional associated to the measure $\dfrac{x^{\alpha+1}}{x-\xi} e^{-x}dx$ and therefore, we have $c^{\alpha+1}(.)=\tilde c_1((x-\xi) .)$. Other notations will be used:
\begin{eqnarray}
k_n^{\alpha}&=&c^{\alpha}((L_n^{\alpha}(x))^2)= n! \Gamma(n+\alpha +1), \label{eq:D1.0} \\
\hat k_n &=& c_1((T_n(x,\xi,M))^2). \nonumber
\end{eqnarray}
Relation (\ref{eq:1.1}) becomes
\begin{equation}\label{eq:D1.1}
T_n(x,\xi,M)= L_n^{\alpha+1}(x) - \sigma_n L_{n-1}^{\alpha+1}(x).
\end{equation}
 Let us prove the following property giving a relation between $T_{n+1}(t,\xi,M)$ and the reproducing kernel $N_n^{\alpha+1}(x,t)$.
%property
 \begin{property}\label{2.12}
 The monic polynomial $T_{n+1}(t,\xi,M)$, orthogonal with respect to $c_1$, satisfies the following relation
 \begin{equation}\label{eq:D1.2}
T_{n+1}(t,\xi,M)= -\frac{k_n^{\alpha+1}}{c_1(L_n^{\alpha+1}(x))} \left(1-(t-\xi)c_1(N_n^{\alpha+1}(x,t))\right)
 \end{equation}
 where the linear functional $c_1$ acts on the variable $x$ and $N_n^{\alpha+1}(x,t)$ is the reproducing kernel 
 $$
 N_n^{\alpha+1}(x,t)= \sum_{j=0}^n \frac{L_j^{\alpha+1}(x)L_j^{\alpha+1}(t)}{k_j^{\alpha+1}}.
 $$
  \end{property}
 \noindent {\bf Proof.}\\
 $T_{n+1}(x,\xi,M)$ is written in the following basis $\{1,\{(x-\xi)L_j^{\alpha+1}(x)\}_{j=0}^n\}$. Thus
 \begin{equation}\label{eq:D2.1}
 T_{n+1}(x,\xi,M)=T_{n+1}(\xi,\xi,M)+ \sum_{j=0}^n \theta_{j,n} (x-\xi)L_j^{\alpha+1}(x), \quad \theta_{n,n}=1.
 \end{equation}
 (\ref{eq:D2.1}) is multiplied by $L_i^{\alpha+1}(x)$ and $c_1$ is applied, for $i=0,\ldots,n$. We obtain
 \begin{eqnarray*}
c_1(L_i^{\alpha+1}(x)  T_{n+1}(x,\xi,M))& = & 0, \quad i=0,\ldots,n,\\
  &=& \tilde c_1(L_i^{\alpha+1}(x) ) T_{n+1}(\xi,\xi,M) + \sum_{j=0}^n \theta_{j,n} c^{\alpha+1}(L_i^{\alpha+1}(x)L_j^{\alpha+1}(x))\\
 &&+ M L_i^{\alpha+1}(\xi)  T_{n+1}(\xi,\xi,M)\\
 &=& T_{n+1}(\xi,\xi,M)\left( \tilde c_1(L_i^{\alpha+1}(x) )+ M L_i^{\alpha+1}(\xi)\right) + \theta_{i,n} k_i^{\alpha+1}\\
 &=&T_{n+1}(\xi,\xi,M) c_1(L_i^{\alpha+1}(x) ) + \theta_{i,n} k_i^{\alpha+1}.
 \end{eqnarray*}
 Hence
 \begin{equation}\label{eq:D2.2}
 \theta_{i,n}=- \frac{T_{n+1}(\xi,\xi,M)}{k_i^{\alpha+1}}c_1(L_i^{\alpha+1}(x) ), \quad i=0,\ldots,n,
  \end{equation}
 and since $\theta_{n,n}=1$ we have
 \begin{equation}\label{eq:D4.2}
T_{n+1}(\xi,\xi,M)=- \frac{k_n^{\alpha+1}}{c_1(L_n^{\alpha+1}(x))}.
 \end{equation}
 Thus
  \begin{eqnarray*}
 T_{n+1}(t,\xi,M)&= &- \frac{k_n^{\alpha+1}}{c_1(L_n^{\alpha+1}(x))} \left(1-(t-\xi) \sum_{j=0}^n \frac{c_1(L_j^{\alpha+1}(x))L_j^{\alpha+1}(t)}{k_j^{\alpha+1}}\right)\\
 &=& - \frac{k_n^{\alpha+1}}{c_1(L_n^{\alpha+1}(x))} \left(1-(t-\xi)c_1(N_n^{\alpha+1}(x,t))\right).
  \end{eqnarray*}
$\square$\\
Now relations giving $\hat k_n$ and $\sigma_n$ can be proved.
%property
\begin{property}\label{2.16}
The square norm $\hat k_n$ of the polynomials $T_n(x,\xi,M)$ and $\sigma_n$ are given by
 \begin{eqnarray}
\hat k_0&=&c_1(1), \nonumber\\
\hat k_n&=& - k_{n-1}^{\alpha+1} \frac{c_1(L_n^{\alpha+1}(x))}{c_1(L_{n-1}^{\alpha+1}(x))}, \quad \forall nÊ\ge 1,\label{eq:D4.3}\\
\sigma_n &=& \frac{c_1(L_n^{\alpha+1}(x))}{c_1(L_{n-1}^{\alpha+1}(x))},  \quad \forall nÊ\ge 1. \label{eq:D4.4}
\end{eqnarray}
\end{property}
 \noindent {\bf Proof.}
 \begin{eqnarray*}
\hat k_n&=& c_1(L_n^{\alpha+1}(x)T_{n}(x,\xi,M) )\\
&=& c_1(L_n^{\alpha+1}(x)T_{n}(\xi,\xi,M) ) +  c_1(L_n^{\alpha+1}(x)\sum_{j=0}^{n-1} \theta_{j,n} (x-\xi)L_j^{\alpha+1}(x) ) \\
&=& c_1(L_n^{\alpha+1}(x)T_{n}(\xi,\xi,M) ).
 \end{eqnarray*}
By using (\ref{eq:D4.2}), (\ref{eq:D4.3}) holds. Now
\begin{eqnarray*}
\hat k_n&=& c_1((x-\xi) L_{n-1}^{\alpha+1}(x)T_{n}(x,\xi,M) )\\
&=& c^{\alpha+1}( L_{n-1}^{\alpha+1}(x)T_{n}(x,\xi,M) )\\
&=& -\sigma_n k_{n-1}^{\alpha+1} \mbox{ by using (\ref{eq:D1.1})}.\\
\sigma_n&=& - \frac{\hat k_n}{k_{n-1}^{\alpha+1}} = \frac{c_1(L_n^{\alpha+1}(x))}{c_1(L_{n-1}^{\alpha+1}(x))} \mbox{ by using (\ref{eq:D4.3}). }
 \end{eqnarray*}
$\square$\\
At last we have the following theorem concerning the Markov-Bernstein inequalities.
%theorem
\begin{theorem}\label{2.17}
The following Markov-Bernstein inequality holds:
$$
c_1((p^\prime)^2) \le \frac{1}{\mu_{1,n}} c_0(p^2), \quad \forall p \in \mathcal{P}_n
$$
where $\mu_{1,n}$ is the smallest zero of the polynomials $A_n(\lambda)$ satisfying the following three term recurrence relation
\begin{equation}\label{eq:D5.1}
A_n (\lambda) = (\lambda +B_n)A_{n-1} (\lambda)-C_n A_{n-2} (\lambda), \quad \forall n\ge 2.
\end{equation}
where 
\begin{eqnarray}
B_n&=& \frac{(n-1)(n+\alpha )}{n} \frac{c_1(L_{n-2}^{\alpha+1}(x))}{c_1(L_{n-1}^{\alpha+1}(x))} + \frac{1}{n-1} \frac{c_1(L_{n-1}^{\alpha+1}(x))}{c_1(L_{n-2}^{\alpha+1}(x))}, \label{eq:D5.2}\\
C_n&=& \frac{(n+\alpha -1)(n -2)}{(n-1)^2}\frac{c_1(L_{n-1}^{\alpha+1}(x))c_1(L_{n-3}^{\alpha+1}(x))}{(c_1(L_{n-2}^{\alpha+1}(x))^2}.\label{eq:D5.3}
\end{eqnarray}
$A_0(\lambda)=1$ and 
 \begin{equation}\label{eq:D5.4}
A_1(\lambda)=\lambda-\frac{ k_1^{ \alpha}}{c_1(1)} = \lambda-\frac{ \Gamma(\alpha +2)}{c_1(1)}.
\end{equation}
\end{theorem}
\noindent {\bf Proof.}\\
By using (\ref{eq:D4.3}), (\ref{eq:D4.4}) and (\ref{eq:D1.0}) we have
\begin{eqnarray*}
B_n&=& - \frac{k_{n}^{\alpha}}{n^2 \hat k_{n-1}}- \frac{\sigma_{n-1}^2 }{({n-1})^2}\frac{ k_{n-1}^{\alpha}}{\hat k_{n-1}}\\
&=&   \frac{(n-1)(n+\alpha )}{n} \frac{c_1(L_{n-2}^{\alpha+1}(x))}{c_1(L_{n-1}^{\alpha+1}(x))} + \frac{1}{n-1} \frac{c_1(L_{n-1}^{\alpha+1}(x))}{c_1(L_{n-2}^{\alpha+1}(x))},\\
C_n &=&  \frac{\sigma_{n-1}^2 (k_{n-1}^{\alpha})^2}{(n-1)^4 \hat k_{n-1} \hat k_{n-2}} \\
&=& \frac{(n+\alpha -1)(n -2)}{(n-1)^2}\frac{c_1(L_{n-1}^{\alpha+1}(x))c_1(L_{n-3}^{\alpha+1}(x))}{(c_1(L_{n-2}^{\alpha+1}(x))^2}.
 \end{eqnarray*}
 $\square$\\
Note that $c_1(L_{i}^{\alpha+1}(x)) = \tilde c_1(L_{i}^{\alpha+1}(x)) + M L_{i}^{\alpha+1}(\xi)$. Therefore to compute 
$c_1(L_{i}^{\alpha+1}(x))$ is reduced to compute $\tilde c_1(L_{i}^{\alpha+1}(x))$.
 %remark
 \begin{remark}\label{2.18}
 If $\xi=0$ and $M=0$, (\ref{eq:D5.1}) is 
 \begin{equation}\label{eq:D6.1}
A_n(\lambda)=(\lambda -2-\frac{\alpha}{n}) A_{n-1}(\lambda)-(1+\frac{\alpha}{n-1})  A_{n-2}(\lambda), \quad \forall n\ge 2,
 \end{equation}
 Indeed $c_1(L_{i}^{\alpha+1}(x)) = \tilde c_1(L_{i}^{\alpha+1}(x))$ and by using the relation $L_{i}^{\alpha}(x) = L_{i}^{\alpha+1}(x) +i L_{i-1}^{\alpha+1}(x)$ we have 
  \begin{equation}\label{eq:D6.2}
  \tilde c_1(L_{i}^{\alpha+1}(x)) = - i  \tilde c_1(L_{i-1}^{\alpha+1}(x)) = (-1)^i  i!\tilde c_1(L_{0}^{\alpha+1}(x)) = (-1)^i  i! \Gamma(\alpha +1).
 \end{equation}
 This last relation, used in (\ref{eq:D5.2}) and (\ref{eq:D5.3}), provides the values of $B_n$ and $C_n$ of (\ref{eq:D6.1}). (\ref{eq:D6.1}) is the relation obtained in the study of Markov-Bernstein inequalities (see Remark \ref{2.3}) for the Laguerre-Sonin measure.
 \end{remark}
  %remark
 \begin{remark}\label{2.19}
 The Laguerre function of the second kind $Q_{n}^{\alpha+1}(\xi)$ can be used to replace the terms $c_1(L_{n}^{\alpha+1}(x))$ in $B_n$ and $C_n$ given in Theorem \ref{2.17}. $Q_{n}^{\alpha}(\xi)$ is defined as follows (see \cite{NU} Relation (5)  page 287) for $\xi <0$:
\begin{eqnarray*}
Q_{n}^{\alpha}(\xi)&=&\frac{1}{e^{-\xi}\xi^\alpha}\int_{0}^\infty \frac{\tilde L_{n}^{\alpha}(x)}{x-\xi} e^{-x}x^\alpha dx\\
&=&e^{-i\pi \alpha}\Gamma(n+\alpha+1) e^\xi G(n+\alpha+1, \alpha+1;-\xi)
\end{eqnarray*}
where $\tilde L_{n}^{\alpha}(x)$ is the non-monic Laguerre polynomial ($\tilde L_{n}^{\alpha}(x)=\frac{(-1)^n}{n!} L_{n}^{\alpha}(x)$) and $G$ is the confluent hypergeometric function of the second kind given by (see \cite{NU} Relation (16)  page 272)
$$
G(\alpha,\gamma;z)= \frac{\Gamma(1-\gamma)}{\Gamma(\alpha+1-\gamma)} F(\alpha,\gamma;z) + \frac{\Gamma(\gamma-1)}{\Gamma(\alpha)} z^{1-\gamma} F(\alpha+1-\gamma,2-\gamma;z)
$$
and 
$F(\alpha,\gamma;z)=\sum_{\nu=0}^\infty \frac{(\alpha)_\nu}{(\gamma)_\nu \nu!} z^\nu$.\\
Unfortunately with this function $Q_{n}^{\alpha+1}(\xi)$ the coefficients $B_n$ and $C_n$ do not become easier to use.
 \end{remark}
  %remark
\begin{remark}\label{2.20}
Now we will use the asymptotic behavior of $L_{n}^{\alpha+1}(\xi)$ and of $Q_{n}^{\alpha+1}(\xi)$ to compute $\sigma$ and $\tau$ given in Blumenthal Theorem \ref{1.5}. The asymptotic behavior of $L_{n}^{\alpha+1}(\xi)$ was given by Perron (see Szeg\"o \cite{S} Theorem 8.22.3 page 199). For $\xi <0$ we have
 \begin{equation}\label{eq:LL2.1}
L_{n}^{\alpha+1}(\xi)= \frac{e^{\xi/2}}{2 \sqrt{\pi}} \frac{e^{2\sqrt{-n\xi}} n^{\alpha/2+1/4}}{(-\xi)^{\alpha/2+3/4}}(1+O(\frac{1}{\sqrt{n}})).
\end{equation}
The asymptotic behavior of $\tilde c_1(L_{n}^{\alpha+1}(x))$ can be found in the paper of Meijer et al (\cite{MPP} Lemma 4.3). For $\xi <0$ we have
 \begin{equation}\label{eq:LL2.2}
\tilde c_1(L_{n}^{\alpha+1}(x)) \sim e^{-\xi/2} \sqrt{\pi} e^{-2\sqrt{-n\xi}} n^{\alpha/2+1/4}(-\xi)^{\alpha/2+1/4}.
\end{equation}
If $M>0$, then $c_1(L_{n}^{\alpha+1}(x))=\tilde c_1(L_{n}^{\alpha+1}(x))+M L_{n}^{\alpha+1}(\xi) \sim M L_{n}^{\alpha+1}(\xi)$ and 
$$
\frac{c_1(L_{n-1}^{\alpha+1}(\xi))}{c_1(L_{n-2}^{\alpha+1}(\xi))}\sim \frac{L_{n-1}^{\alpha+1}(\xi)}{L_{n-2}^{\alpha+1}(\xi)}\sim (\frac{n-1}{n-2})^{\alpha/2+1/4}  e^{\sqrt{-\xi/(n-1)}}\to 1 \mbox{ when } n \to \infty.
$$
Therefore $\lim_{n \to \infty} B_n = \infty$ and $\lim_{n \to \infty} C_n =1$. Hence $\sigma=\tau=\infty$.\\
If $M=0$, then $c_1(L_{n}^{\alpha+1}(x))=\tilde c_1(L_{n}^{\alpha+1}(x))$ and 
$$
\frac{c_1(L_{n-1}^{\alpha+1}(\xi))}{c_1(L_{n-2}^{\alpha+1}(\xi))}\sim \frac{L_{n-1}^{\alpha+1}(\xi)}{L_{n-2}^{\alpha+1}(\xi)}\sim (\frac{n-1}{n-2})^{\alpha/2+1/4}  e^{-\sqrt{-\xi/(n-1)}}\to 1 \mbox{ when } n \to \infty.
$$
Therefore we have the same result.\\
 \end{remark}
%property
\begin{property}\label{2.21}
If $\xi=0$, then
$$
\lim_{n \to \infty} \mu_{1,n}=0.
$$
\end{property}
\noindent {\bf Proof.}\\
In this case $\tilde c_1=c_0=c^\alpha$. By using the relation $L_{j}^{\alpha}(x)=L_{j}^{\alpha+1}(x)+jL_{j-1}^{\alpha+1}(x)$
for $j=n,n-1,\ldots,1$ we obtain
$c^\alpha (L_{n}^{\alpha+1}(x))=(-1)^n n! k_0^\alpha=(-1)^n n!  \Gamma(\alpha+1)$.\\
For $M>0$ we have $c_1 (L_{n}^{\alpha+1}(x))= (-1)^n (n!  \Gamma(\alpha+1)+M (\alpha+1)_n)$. Therefore
\begin{eqnarray*}
c_1 (L_{n}^{\alpha+1}(x))&\sim&(-1)^n n!  \Gamma(\alpha+1) \quad \mbox{ if } -1<\alpha<0,\\
&=&(-1)^n n!  (M +1) \quad \mbox{ if } \alpha=0,\\
&\sim&(-1)^n M (\alpha+1)_n \quad \mbox{ if } \alpha>0.
\end{eqnarray*}
In the three cases we have $\lim_{n \to \infty} -B_n = 2$ and $\lim_{n \to \infty} C_n =1$. 
We use Blumenthal Theorem \ref{1.5}: we have $\sigma=0$ and $\tau =4$. Hence the property holds by using Remark \ref{1.6}.\\
If $M=0$, the property also holds (see Remark \ref{2.18}). In this case the values of $\sigma=0$ and $\tau =4$ are obvious.
$\square$

%**************************************************************************************************************************	
%_______________________________________________section 3
%_________________________________________________________________
%**************************************************************************************************************************

\section{Jacobi case}

%**************************************************************************************************************************

%_____________________________________________section 3.1__________________________

%**************************************************************************************************************************
\subsection{Jacobi case (a)} 

The measure associated to $c_0$ is $\mid x-\xi \mid (1-x)^{\alpha-1}(1+x)^{\beta-1}dx$ with $\alpha >0$, $\beta >0$ and $\mid \xi \mid > 1$, and the one associated to $c_1$ is the Jacobi measure $(1-x)^{\alpha}(1+x)^{\beta}dx$. Therefore $T_n(x)$ is the monic Jacobi polynomial $ P_n^{(\alpha,\beta)}(x)$. For more convenience the measure associated to $c_0$ will be written as $\varepsilon ( x-\xi ) (1-x)^{\alpha-1}(1+x)^{\beta-1}dx$ with $\varepsilon=1$ if $\xi<-1$ and $\varepsilon=-1$ if $\xi>1$.\\
The monic polynomials, orthogonal with respect to $c_0$, will be denoted by $ P_n(x,\xi)$. Their expression, depending on the monic Jacobi polynomials $P_n^{(\alpha-1,\beta-1)}(x)$, has already also been given in the thesis of P\'erez \cite{P} (Corollary 5.3.3 and Theorem 5.3.4). Also here we give a direct proof of this expression because we will also obtain the square norm $k_n^{(0)}$ and $\sigma_n$.
%theorem
\begin{theorem}\label{3.1}
The monic polynomials $P_n(x,\xi)$, orthogonal with respect to $c_0$, satisfy the following relation:
\begin{equation}\label{eq:17.1}
P_n(x,\xi)= N_n^{(\alpha-1,\beta-1)}(x,\xi) \frac{k_n^{(\alpha-1,\beta-1)}}{P_n^{(\alpha-1,\beta-1)}(\xi)}
\end{equation}
where $N_n^{(\alpha-1,\beta-1)}(x,\xi)$ is the reproducing kernel associated to the measure $(1-x)^{\alpha-1}(1+x)^{\beta-1}dx$.
\begin{equation*}
%\label{eq:17.2}
N_n^{(\alpha-1,\beta-1)}(x,\xi) =\sum_{j=0}^n \frac{P_j^{(\alpha-1,\beta-1)}(x)P_j^{(\alpha-1,\beta-1)}(\xi)}{k_j^{(\alpha-1,\beta-1)}}. 
\end{equation*}
The square norm $k_n^{(0)}$ of $P_n(x,\xi)$ and $\sigma_n$ are given by
\begin{eqnarray}
k_n^{(0)}&=&-\varepsilon  k_n^{(\alpha-1,\beta-1)} \frac{P_{n+1}^{(\alpha-1,\beta-1)}(\xi)}{P_{n}^{(\alpha-1,\beta-1)}(\xi)}, \quad \forall n \ge 0, \mbox{ with } \varepsilon =1 \mbox{ if } \xi <-1  \mbox{ and } \varepsilon =-1 \mbox{ if } \xi >1 \label{eq:18.1}\\
\sigma_n&=&\frac{4n(n+\alpha)(n+\beta)(n+\alpha+\beta-1)}{(2n+\alpha+\beta-1)(2n+\alpha+\beta)^2(2n+\alpha+\beta+1)}\frac{P_{n}^{(\alpha-1,\beta-1)}(\xi)}{P_{n+1}^{(\alpha-1,\beta-1)}(\xi)}, \quad \forall n \ge 1 \label{eq:18.2}
\end{eqnarray}
where $k_n^{(\alpha-1,\beta-1)}$ is the square norm of the polynomials $P_n^{(\alpha-1,\beta-1)}(x)$ with respect to $(1-x)^{\alpha-1}(1+x)^{\beta-1}dx$.
\end{theorem}
\noindent {\bf Proof.}\\
$(x-\xi)P_n(x,\xi)$ is written in the basis of monic Jacobi polynomials $P_j^{(\alpha-1,\beta-1)}(x)$.
\begin{equation}\label{eq:18.3}
(x-\xi)P_n(x,\xi)= \sum_{j=0}^{n+1} \theta_{j,n+1} P_j^{(\alpha-1,\beta-1)}(x), \quad \theta_{n+1,n+1}=1.
\end{equation}
Let us denote by $c^{(\alpha-1,\beta-1)}$ the linear functional defined from the Jacobi measure $(1-x)^{\alpha-1}(1+x)^{\beta-1}dx$ on the support $\Omega=]-1,1[$. We have
\begin{eqnarray*}
c^{(\alpha-1,\beta-1)}(\varepsilon P_i^{(\alpha-1,\beta-1)}(x) (x-\xi)P_n(x,\xi))&=&c_0(P_i^{(\alpha-1,\beta-1)}(x) P_n(x,\xi)) =0, \quad i=0,\ldots,n-1\\
&=&\varepsilon  \sum_{j=0}^{i} \theta_{j,n+1} c^{(\alpha-1,\beta-1)}(P_i^{(\alpha-1,\beta-1)}(x) P_j^{(\alpha-1,\beta-1)}(x))\\
&=& \theta_{i,n+1} k_i^{(\alpha-1,\beta-1)}.
\end{eqnarray*}
Therefore $\theta_{i,n+1}=0$ for $i=0,\ldots,n-1$ and 
\begin{equation}\label{eq:18.4}
\theta_{n,n+1}= \frac{k_n^{(0)}}{\varepsilon k_n^{(\alpha-1,\beta-1)}}.
\end{equation}
We put $x=\xi$ in (\ref{eq:18.3}). Then we obtain
\begin{equation}\label{eq:18.5}
\theta_{n,n+1}= - \frac{P_{n+1}^{(\alpha-1,\beta-1)}(\xi)}{P_{n}^{(\alpha-1,\beta-1)}(\xi)}.
\end{equation}
The expression (\ref{eq:18.1}) of the square norm $k_n^{(0)}$ is deduced from (\ref{eq:18.4}) and (\ref{eq:18.5}).\\
From (\ref{eq:18.3}) we have
\begin{equation}\label{eq:19.1}
P_n(x,\xi)= \left( P_{n+1}^{(\alpha-1,\beta-1)}(x)- \frac{P_{n+1}^{(\alpha-1,\beta-1)}(\xi)}{P_{n}^{(\alpha-1,\beta-1)}(\xi)} P_{n}^{(\alpha-1,\beta-1)}(x) \right) \frac{1}{x-\xi}.
\end{equation}
Again, by using the definition of the reproducing kernel $N_n^{(\alpha-1,\beta-1)}(x,\xi)$ associated to the Christoffel-Darboux formula, 
(\ref{eq:17.1}) is obtained.\\
Note that (\ref{eq:19.1}) is also given in Lemma 3.2 of the paper of Pan \cite{PAN}.

Now let us determine explicitly the sequence $\{\sigma_n \}_{n\ge1}$ as a function depending on $\xi$. From (\ref{eq:1.1}), (\ref{eq:7.5}) and (\ref{eq:17.1}) we have
\begin{eqnarray}
\sigma_n &=& \frac{n}{n+1} \frac{P_{n+1}^\prime(x,\xi)-(n+1)P_n^{(\alpha,\beta)}(x)}{P_n^\prime(x,\xi)} \nonumber\\
&=&  \frac{n}{(n+1)P_n^\prime(x,\xi)} \sum_{j=1}^{n} j P_{j-1}^{(\alpha,\beta)}(x) \frac{P_{j}^{(\alpha-1,\beta-1)}(\xi)}{P_{n+1}^{(\alpha-1,\beta-1)} (\xi)} \frac{k_{n+1}^{(\alpha-1,\beta-1)}}{k_j^{(\alpha-1,\beta-1)}}. \label{eq:19.2}
\end{eqnarray}
(\ref{eq:18.2}) is deduced by using the ratio of the coefficients of $x^{n+1}$ in the numerator and the denominator of (\ref{eq:19.2}).
\begin{equation}\label{eq:19a.1}
\sigma_n=\frac{n}{(n+1)} \frac{P_{n}^{(\alpha-1,\beta-1)}(\xi)}{P_{n+1}^{(\alpha-1,\beta-1)}(\xi)}\frac{k_{n+1}^{(\alpha-1,\beta-1)}}{k_n^{(\alpha-1,\beta-1)}}.
\end{equation}
$\square$\\
Note that (\ref{eq:18.2}) is also given in \cite{MP}.\\
Moreover note that (\ref{eq:19a.1}) and (\ref{eq:18.1}) give a relation between $k_n^{(0)}$ and $\sigma_n$.
\begin{equation}\label{eq:19a.2}
\sigma_n=- \varepsilon \frac{n}{(n+1)} \frac{k_{n+1}^{(\alpha-1,\beta-1)}}{k_n^{(0)}}.
\end{equation}
At last we have the following theorem concerning the Markov-Bernstein inequalities.
%theorem
\begin{theorem}\label{3.2}
The following Markov-Bernstein inequality holds:
$$
c_1((p^\prime)^2) \le \frac{1}{\mu_{1,n}} c_0(p^2), \quad \forall p \in \mathcal{P}_n
$$
where $\mu_{1,n}$ is the smallest zero of the polynomials $A_n(\lambda)$ satisfying the following three term recurrence relation
\begin{equation}\label{eq:20.1}
A_n (\lambda) = (\lambda +B_n)A_{n-1} (\lambda)-C_n A_{n-2} (\lambda), \quad \forall n\ge 2.
\end{equation}
where 
\begin{eqnarray}
B_n&=& \frac{\varepsilon}{n(n+\alpha +\beta-1)} \left( \xi-\frac{(\beta-\alpha)(\alpha +\beta-2)}{(2n+\alpha +\beta-2)(2n+\alpha +\beta)}   \right),\label{eq:20.2}\\
C_n&=& \frac{4(n+\alpha -1)(n +\beta-1)}{(n-1)(n+\alpha +\beta-1)(2n+\alpha +\beta-1)(2n+\alpha +\beta-2)^2(2n+\alpha +\beta-3)}.\label{eq:20.3}
\end{eqnarray}
$A_0(\lambda)=1$ and 
\begin{eqnarray*}
A_1(\lambda)&= &\lambda+\frac{ \varepsilon}{\alpha +\beta}  \frac{P_{2}^{(\alpha-1,\beta-1)}(\xi)}{P_{1}^{(\alpha-1,\beta-1)}(\xi)}\\
&=& \lambda+\frac{ \varepsilon}{\alpha +\beta} \left(\xi- \frac{(\beta-\alpha)(\alpha +\beta-2)}{(2+\alpha +\beta)(\alpha +\beta)}-\frac{4 \alpha \beta}{(\alpha +\beta)_2(\xi(\alpha +\beta)-\beta +\alpha )}  \right) \\
&=& \lambda+\varepsilon \frac{(1+\alpha +\beta)_2 \xi^2+ 2 (\alpha -\beta)(1+\alpha +\beta)\xi + (\alpha -\beta)^2 -(2+\alpha +\beta)}{(1+\alpha +\beta)_2(\xi(\alpha +\beta)-\beta +\alpha )} 
\end{eqnarray*}
\end{theorem}
\noindent {\bf Proof.}\\
By using (\ref{eq:19a.2}),
$$
C_n= \frac{\sigma_{n-1}^2 (k_{n-1}^{(0)})^2}{(n-1)^4 k_{n-1}^{(\alpha ,\beta)}k_{n-2}^{(\alpha ,\beta)}}  =  \frac{(k_n^{(\alpha-1,\beta-1)})^2}{n^2(n-1)^2 k_{n-1}^{(\alpha ,\beta)}k_{n-2}^{(\alpha ,\beta)}}.
$$
Hence (\ref{eq:20.3}) holds.
\begin{eqnarray*}
B_n&=&- \frac{k_{n}^{(0)}}{n^2 k_{n-1}^{(\alpha ,\beta)}}- \frac{\sigma_{n-1}^2 }{({n-1})^2}\frac{ k_{n-1}^{(0)}}{k_{n-1}^{(\alpha ,\beta)}}\\
&=& \frac{ \varepsilon}{n^2} \frac{k_n^{(\alpha-1,\beta-1)}}{ k_{n-1}^{(\alpha ,\beta)}} \left( \frac{P_{n+1}^{(\alpha-1,\beta-1)}(\xi)}{P_{n}^{(\alpha-1,\beta-1)}(\xi)} + \frac{P_{n-1}^{(\alpha-1,\beta-1)}(\xi)}{P_{n}^{(\alpha-1,\beta-1)}(\xi)} \frac{k_n^{(\alpha-1,\beta-1)}}{ k_{n-1}^{(\alpha-1,\beta-1)}} \right)
\end{eqnarray*}
by using (\ref{eq:19a.1}) and (\ref{eq:18.1}).\\
Now, by using the three term recurrence relation (\ref{eq:7.6}) the quantity inside the brackets is equal to
$$
\xi-\frac{(\beta-\alpha)(\alpha +\beta-2)}{(2n+\alpha +\beta-2)(2n+\alpha +\beta)}.
$$
Hence (\ref{eq:20.2}) holds.
At last 
\begin{eqnarray*}
A_1(\lambda)&=& \lambda-\frac{k_{1}^{(0)}}{k_{0}^{(\alpha ,\beta)}}\\
&=&\lambda+\frac{ \varepsilon}{\alpha +\beta}  \frac{P_{2}^{(\alpha-1,\beta-1)}(\xi)}{P_{1}^{(\alpha-1,\beta-1)}(\xi)}\\
&=& \lambda+\frac{ \varepsilon}{\alpha +\beta} \left(\xi- \frac{(\beta-\alpha)(\alpha +\beta-2)}{(2+\alpha +\beta)(\alpha +\beta)}-\frac{4 \alpha \beta}{(1+\alpha +\beta)(\alpha +\beta)(\xi(\alpha +\beta)-\beta +\alpha )}  \right)
\end{eqnarray*}
by using the three term recurrence relation (\ref{eq:7.1}) satisfied by the polynomials $P_{2}^{(\alpha-1,\beta-1)}(\xi)$, $P_{1}^{(\alpha-1,\beta-1)}(\xi)$ and $P_{0}^{(\alpha-1,\beta-1)}(\xi)$.
$\square$
%property
\begin{property}\label{3.2a}
$$
\lim_{n \to \infty} \mu_{1,n}=0.
$$
\end{property}
\noindent {\bf Proof.}\\
We use Blumenthal Theorem \ref{1.5}. $\lim_{n \to \infty} B_n = 0$ and $\lim_{n \to \infty} C_n =0$. Therefore $\sigma=0$ and $\tau =0$. Hence the property holds by using Remark \ref{1.6}.
$\square$

%**************************************************************************************************************************

%_____________________________________________section 3.2_______________________________

%**************************************************************************************************************************
\subsection{Jacobi case (b) and (c)} 
 
In the case (b) the measure associated to $c_0$ is $(1+x)^{\beta-1}dx + M \delta(1)$ with $\beta >0$ and $M\ge 0$, and the one associated to $c_1$ is the Jacobi measure $(1+x)^{\beta}dx$. In the case (c) the measure associated to $c_0$ is $ (1-x)^{\alpha-1}dx+  M \delta(-1)$ with $\alpha >0$ and $M\ge 0$, and the one associated to $c_1$ is the Jacobi measure $(1-x)^{\alpha}dx$. For more convenience and to give unique proofs for both cases (b) and (c) we will adopt a specific notation. The measure associated to $c_0$ will be $(1+\varepsilon x)^{\gamma-1}dx + M \delta(\varepsilon)$ with $\gamma >0$ and $M\ge 0$, and the one associated to $c_1$ will be $ (1+\varepsilon x)^{\gamma}dx$. If $\varepsilon=1$ then $\gamma=\beta$, and if $\varepsilon=-1$ then $\gamma=\alpha$. The monic polynomials, orthogonal with respect to $c_0$, will be denoted by $P_n(x,M)$. The monic polynomials $T_n(x)$, orthogonal with respect to $c_1$, will be denoted by $P_n^{(\gamma)}$. If $\varepsilon=1$, this polynomial is the monic Jacobi polynomial $P_n^{(0,\beta)}(x)$, and if $\varepsilon=-1$,  this polynomial is the monic Jacobi polynomial $P_n^{(\alpha,0)}(x)$.
 Thanks to this specific notation we have unique expressions for the three term recurrence relation satisfied by the polynomials $P_n^{(\gamma)}$, for the square norm $k_n^{(\gamma)}$, and for the explicit formula of $P_n^{(\gamma)}$ corresponding to (\ref{eq:7.3}) and (\ref{eq:7.4}). 
\begin{eqnarray} 
P_{n+1}^{(\gamma)}(x)&=& (x-\frac{\varepsilon \gamma^2}{(2n+\gamma)(2n+\gamma+2)})P_n^{(\gamma)}(x) \nonumber\\
&&- \frac{4 n^2 (n+\gamma)^2}{(2n+\gamma)^2(2n+\gamma+1)(2n+\gamma-1)} P_{n-1}^{(\gamma)}(x), \label{eq:23.1}\\
k_{n}^{(\gamma)}&=&2^{2n+\gamma+1} (n!)^2 \frac{(\Gamma(n+\gamma +1))^2}{(2n+\gamma +1)(\Gamma(2n+\gamma +1))^2}, \label{eq:23.1a}\\
P_n^{(\gamma)}&=& \frac{\varepsilon^n 2^n n!}{\Gamma(\gamma +2n+1)}\sum_{m=0}^n \left(\begin{array}{c}n\\m\end{array}\right)\frac{\Gamma(\gamma +n+m+1)}{m!}(\frac{\varepsilon x-1}{2})^m. \label{eq:23.2}
\end{eqnarray}
Now $P_n(x,M)$ is written in the basis of monic polynomials $P_j^{(\gamma-1)}(x)$ orthogonal with respect to $\hat c_0$ which is the linear functional associated to the measure  $(1+\varepsilon x)^{\gamma-1}dx$ on the support $]-1,1[$. 
$$
P_n(x,M) = \sum_{j=0}^n \theta_{j,n} P_j^{(\gamma-1)}(x), \quad \theta_{n,n}=1.
$$
%property
\begin{property}\label{3.7}
The explicit expressions of $P_n(x,M)$ in the basis of monic Jacobi polynomials $P_j^{(\gamma-1)}(x)$, of the square norm $k_{n}^{(0)}$ of  $P_n(x,M)$, and of $\sigma_n$ are given by
\begin{eqnarray}
P_n(x,M)&=& P_n^{(\gamma-1)}(x)-  \frac{M n! 2^n \varepsilon^n}{2^\gamma+n M (n-1+ \gamma)} 
\frac{\Gamma(n+\gamma)}{\Gamma(2n+\gamma)}  \sum_{j=0}^{n-1}  \frac{ \varepsilon^j }{2^j}\frac{(\gamma +j)_{j+1}}{j!} P_j^{(\gamma-1)}(x). \label{eq:J.1}\\\
k_{n}^{(0)} & =&  \left(\frac{\Gamma(n+\gamma)}{\Gamma(2n+\gamma)}\right)^2 \frac{2^{2n+\gamma}(n!)^2} {2n+\gamma}  \frac{2^\gamma+(n+1) M (n+ \gamma)}{2^\gamma+n M (n-1+ \gamma)},  \label{eq:30.1}\\
 \sigma_n& =&    \frac{\varepsilon n} {2n+\gamma+1} \left(  \frac{\gamma} {2n+ \gamma}- \frac{2 M(\gamma+n)}{2^\gamma+(n+1) M (n+ \gamma)} \right) \label{eq:30.2}\\
 &=& \frac{\varepsilon n(\gamma 2^\gamma+M (n+ \gamma)(4n-\gamma(n-1)) )} {(2n+\gamma+1)(2n+ \gamma)(2^\gamma+(n+1) M (n+ \gamma))}.  \label{eq:30.3}
 \end{eqnarray}
\end{property}
\noindent {\bf Proof.}
We use the following orthogonality relations of $P_n(x,M)$
\begin{eqnarray*}
c_0(P_i^{(\gamma-1)}(x) P_n(x,M) )&=& 0 \quad \mbox{ for } i=0,\ldots,n-1 \\
&=& \sum_{j=0}^n \theta_{j,n} \left( \hat c_0(P_i^{(\gamma-1)}(x) P_j^{(\gamma-1)}(x) )+M P_i^{(\gamma-1)}(\xi) P_j^{(\gamma-1)}(\xi) \right)\\
&=&  \frac{\theta_{i,n} 2^{2i+\gamma} (i!)^2 (\Gamma(i+\gamma ))^2}{(2i+\gamma )(\Gamma(2i+\gamma))^2}
+ M \varepsilon ^i 2^i i!  \frac{\Gamma(i+\gamma )}{\Gamma(2i+\gamma)}  \sum_{j=0}^n \theta_{j,n} \varepsilon^j 2^j j!  \frac{(\Gamma(j+\gamma )}{\Gamma(2j+\gamma)} 
\end{eqnarray*}
Thus 
$$
\theta_{i,n} = -\frac{M \varepsilon ^i (2i+\gamma )\Gamma(2i+\gamma)}{ 2^{i+\gamma} i!  \Gamma(i+\gamma )}  \sum_{j=0}^n \theta_{j,n} \varepsilon^j 2^j j!  \frac{(\Gamma(j+\gamma )}{\Gamma(2j+\gamma)}.
$$
Let us set $K= \sum_{j=0}^n \theta_{j,n} \varepsilon^j 2^j j!  \frac{(\Gamma(j+\gamma )}{\Gamma(2j+\gamma)}$.
$$
K= -M \sum_{j=0}^{n-1} (2j+\gamma) 2^{-\gamma} K + \varepsilon^n 2^n n!  \frac{(\Gamma(n+\gamma )}{\Gamma(2n+\gamma)}
$$
by using the expressions of the $\theta_{j,n}$'s in the definition of $K$. Hence
$$
K= \frac{\varepsilon^n 2^n n!}{1+M 2^{-\gamma} n (\gamma+n-1)}  \frac{(\Gamma(n+\gamma )}{\Gamma(2n+\gamma)}
$$
and 
\begin{eqnarray} 
\theta_{j,n} &=& -\varepsilon^{n+j} M \frac{2^{n-j} n!}{j!} \frac{(\gamma+j)_{j+1}}{2^\gamma+n M (n-1+ \gamma)} \frac{\Gamma(n+\gamma)}{\Gamma(2n+\gamma)}, \quad j=0,\ldots,n-1, \label{eq:27.1} \\
\theta_{n,n}&=&1. \nonumber
\end{eqnarray}
Hence Relation (\ref{eq:J.1}) holds.\\
The square norm $k_{n}^{(0)}$ of  $P_n(x,M)$ is given by
\begin{eqnarray*}
k_{n}^{(0)} &=& c_0(( P_n(x,M))^2)\\
&=& \hat c_0(( P_n(x,M))^2) + M ( P_n(\varepsilon,M))^2\\
&=& \sum_{j=0}^n  \theta_{j,n} ^2 k_j^{(\gamma-1)} + M(\sum_{j=0}^n  \theta_{j,n} P_j^{(\gamma-1)}(\varepsilon))^2\\
&=& \left(\frac{\Gamma(n+\gamma)}{\Gamma(2n+\gamma)}\right)^2 \frac{2^{2n+\gamma}(n!)^2} {2n+\gamma}  \frac{2^\gamma+(n+1) M (n+ \gamma)}{2^\gamma+n M (n-1+ \gamma)}
\end{eqnarray*}
by using the expressions of the $\theta_{j,n}$'s, of $k_j^{(\gamma-1)}$ and the value of $P_j^{(\gamma-1)}(\varepsilon)$. Hence (\ref{eq:30.1}) holds.

$\sigma_n$ will be obtained from
$$
\frac{n}{n+1} \frac{P_{n+1}^\prime(x,M)-(n+1)P_n^{(\gamma)}(x)}{P_n^\prime(x,M)}
$$
with the ratio of the coefficients of $x^{n-1}$ in the numerator and the denominator.\\
The coefficient of $x^{n-1}$ in $P_{n +1}^\prime(x,M)$ is
$$
\varepsilon n(n+1) \left( -1+ \frac{2n+1}{2n+\gamma+1} -\frac{2M}{2^\gamma+(n+1) M (n+ \gamma)} \frac{\gamma+n}{\gamma+2n+1}\right).
$$
The coefficient of $x^{n-1}$ in $(n+1)P_n^{(\gamma)}(x)$ is
$$
\varepsilon n(n+1) \left( -1+ \frac{2n}{2n+\gamma} \right).
$$
Therefore (\ref{eq:30.2}) holds.
$\square$\\ 
Like in the Laguerre case (b) a much more complicated way consists to use the following orthogonality relations
\begin{equation}\label{eq:23.3}
c_0((1-\varepsilon x)^i P_n(x,M))= 0, \quad i=0,\ldots,n-1.
\end{equation}
The main interest of this kind of proof is to give explicitly $\hat c_0((1-\varepsilon x)^i P_j^{(\gamma-1)})$ for $i\ge j$, and as a consequence a new identity implying a sum of products of binomial coefficients and Pochhammer coefficients. The expression of $\hat c_0((1-\varepsilon x)^i P_j^{(\gamma-1)})$ for $i\ge j$ is given in the following Lemma. A Corollary will give the new identity.
%lemma
\begin{lemma}\label{3.3}
For $i\ge j$, $\forall j \in \mathbb{N}$,
\begin{equation}\label{eq:23.4}
\hat c_0((1-\varepsilon x)^i P_j^{(\gamma-1)}(x))= \frac{(-\varepsilon)^j 2^{i+j+\gamma} i! j!}{(\gamma +j)_{i+1}}\frac{\Gamma(j+\gamma)}{\Gamma(2j+\gamma)}\left(\begin{array}{c}i\\j\end{array}\right).
\end{equation}
\end{lemma}
\noindent {\bf Proof.}\\
Note that if $i=j$, the right hand side of (\ref{eq:23.4}) exactly is $(-\varepsilon)^j  k_j^{(\gamma-1)}$. Moreover if $j=0$, the right hand side of (\ref{eq:23.4}) is equal to $i! \frac {2^{\gamma+i}}{(\gamma)_{i+1}}$. It exactly is the value of $\hat c_0((1-\varepsilon x)^i)$ which can be proved by recurrence very simply.\\
First of all we write the three term recurrence relation satisfied by the polynomials $ P_j^{(\gamma-1)}$ as
\begin{eqnarray} \label{eq:24.1}
P_{j+1}^{(\gamma-1)}(x)&=& -\varepsilon(1-\varepsilon x-\frac{(2j+\gamma-1)(2j+\gamma+1)- (\gamma-1)^2}{(2j+\gamma-1)(2j+\gamma+1)})P_j^{(\gamma-1)}(x) \nonumber\\
&&- \frac{4 j^2 (j+\gamma-1)^2}{(2j+\gamma-1)^2(2j+\gamma)(2j+\gamma-2)} P_{j-1}^{(\gamma-1)}(x). 
\end{eqnarray}
We begin to prove that Lemma \ref{3.3} holds for $i=j+1$. (\ref{eq:24.1}) is multiplied by $(1-\varepsilon x)^{j}$ and $\hat c_0$ is applied to the result. Since Lemma \ref{3.3} holds for $i=1$ and $j=0$, we assume that it holds for $j=0,\ldots,r-1$ and $i=j+1$, and we prove that it holds for $j=r$ and $i=r+1$. We have
\begin{eqnarray*} 
\hat c_0((1-\varepsilon x)^{r+1} P_{r}^{(\gamma-1)}(x))&=&\frac{(2r+\gamma-1)(2r+\gamma+1)- (\gamma-1)^2}{(2r+\gamma-1)(2r+\gamma+1)} (-\varepsilon)^rk_r^{(\gamma-1)} \nonumber\\
&&- \frac{\varepsilon 4 r^2 (r+\gamma-1)^2}{(2r+\gamma-1)^2(2r+\gamma)(2r+\gamma-2)} \frac{(-\varepsilon)^{r-1} 2^{2r-1+\gamma}  r!(r-1)!}{(\gamma +r-1)_{r+1}}\times\\
&&\frac{\Gamma(r-1+\gamma)}{\Gamma(2r-2+\gamma)}\left(\begin{array}{c}r\\r-1\end{array}\right)\\
&=& \frac{(-\varepsilon)^r 2^{2r+1+\gamma} (r+1)! r!}{(\gamma +r)_{r+1}}\frac{\Gamma(r+\gamma)}{\Gamma(2r+\gamma)}\left(\begin{array}{c}r+1\\r\end{array}\right).
\end{eqnarray*}
Hence Lemma \ref{3.3} holds for $i=r+1$ and $j=r$.

Now we assume that Lemma \ref{3.3} holds $\forall j \in \mathbb{N}$ and $i=j,\ldots,j+r-1$. We prove it for $i=j+r$, $\forall j \in \mathbb{N}$. (\ref{eq:24.1}) is multiplied by $(1-\varepsilon x)^{j+r-1}$ and $\hat c_0$ is applied to the result. $\hat c_0((1-\varepsilon x)^{r+j} P_{j}^{(\gamma-1)}(x))$ is expressed as a function of $\hat c_0((1-\varepsilon x)^{r+j-1} P_{j+1}^{(\gamma-1)}(x))$, $\hat c_0((1-\varepsilon x)^{r+j-1} P_{j}^{(\gamma-1)}(x))$ and $\hat c_0((1-\varepsilon x)^{r+j-1} P_{j-1}^{(\gamma-1)}(x))$. By using the assumption of recurrence, we obtain 
$$
\hat c_0((1-\varepsilon x)^{r+j} P_{j}^{(\gamma-1)}(x))= \frac{(-\varepsilon)^j 2^{2r+j+\gamma} (r+j)! j!}{(\gamma +j)_{r+f+1}}\frac{\Gamma(j+\gamma)}{\Gamma(2j+\gamma)}\left(\begin{array}{c}r+j\\j\end{array}\right).
$$
Hence Lemma \ref{3.3} holds.
$\square$
%corollary
\begin{corollary}\label{3.4}
For $i\ge j$,
\begin{equation}\label{eq:26.1}
 \sum_{m=0}^j (-1)^m  \left(\begin{array}{c}j\\m\end{array}\right) \left(\begin{array}{c}i+m\\i\end{array}\right) (\gamma+i+m+1)_{j-m}(\gamma+j)_m=(-1)^j  \left(\begin{array}{c}i\\j\end{array}\right)  (\gamma)_j.
 \end{equation}
\end{corollary}
\noindent {\bf Proof.}\\
For $i\ge j$, by using (\ref{eq:23.2}) we have
\begin{eqnarray*} 
\hat c_0((1-\varepsilon x)^i P_j^{(\gamma-1)})&=& \frac{\varepsilon^j 2^{j} j!}{\Gamma(2j+\gamma)} \sum_{m=0}^j (-1)^m\left(\begin{array}{c}j\\m\end{array}\right)
\frac{\Gamma(j+m+\gamma)}{2^m m!}  \hat c_0((1-\varepsilon x)^{i+m})\\
&=&  \frac{\varepsilon^j 2^{i+j+\gamma}}{(\gamma)_{i+j+1}}\frac{\Gamma(j+\gamma)}{\Gamma(2j+\gamma)} \sum_{m=0}^j (-1)^m\left(\begin{array}{c}j\\m\end{array}\right) \left(\begin{array}{c}i+m\\i\end{array}\right) \times \\
&& (\gamma+i+m+1)_{j-m}(\gamma+j)_m.
\end{eqnarray*}
by using (\ref{eq:23.4}) for $j=0$. At last Lemma \ref{3.3} implies that (\ref{eq:26.1}) holds.
$\square$\\
From Lemma \ref{3.3} it remains to compute the $\theta_j$'s. This part also is very technical. Since we already have this result, it is not necessary to give it.\\
Finally we have the following theorem concerning the Markov-Bernstein inequalities.
%theorem
\begin{theorem}\label{3.8}
The following Markov-Bernstein inequality holds:
$$
c_1((p^\prime)^2) \le \frac{1}{\mu_{1,n}} c_0(p^2), \quad \forall p \in \mathcal{P}_n
$$
where $\mu_{1,n}$ is the smallest zero of the polynomials $A_n(\lambda)$ satisfying the following three term recurrence relation
\begin{equation}\label{eq:33.1}
A_n (\lambda) = (\lambda +B_n)A_{n-1} (\lambda)-C_n A_{n-2} (\lambda), \quad \forall n\ge 2.
\end{equation}
where 
\begin{eqnarray}
B_n&=&- \frac{2^\gamma+(n+1) M (n+ \gamma)}{2^\gamma+n M (n-1+ \gamma)} \frac{2}{(2n+\gamma)(2n+\gamma-1)} \nonumber\\
&&-
\frac{(\gamma 2^\gamma+M (n-1+ \gamma)(4(n-1)-\gamma(n-2)) )^2}{2(2^\gamma+n M (n-1+ \gamma))(2^\gamma+(n-1) M (n-2+ \gamma))} \times \nonumber \\
&&\frac{1}{(2n+\gamma-1)(2n+\gamma-2)(n-1+\gamma)^2}   ,\label{eq:33.2}\\
C_n&=& \frac{(\gamma 2^\gamma+M (n-1+ \gamma)(4(n-1)-\gamma(n-2)) )^2}{(2^\gamma+(n-1) M (n-2+ \gamma))^2(n+\gamma-1)^2(2n+\gamma-1)(2n+\gamma-2)^2(2n+\gamma-3)}.\label{eq:33.3}
\end{eqnarray}
$A_0(\lambda)=1$ and 
\begin{equation}\label{eq:33.4}
A_1(\lambda)= \lambda-\frac{ 2 (\gamma+1)^2}{\gamma+2}  \frac{2^\gamma+2M(1+\gamma)}{2^\gamma+M\gamma}.
\end{equation}
\end{theorem}
\noindent {\bf Proof.}\\
All these relations are obtained by replacing $k_{n-1}^{(0)}$, $k_{n}^{(0)}$, $k_{n-1}^{(1)}$, $k_{n-2}^{(1)}$ and $\sigma_{n-1}$ in (\ref{eq:6.1}) by their expressions given in the relations  (\ref{eq:30.1}),  (\ref{eq:23.1a}) and  (\ref{eq:30.3}).\\
$\square$
%property
\begin{property}\label{3.8a}
$$
\lim_{n \to \infty} \mu_{1,n}=0.
$$
\end{property}
\noindent {\bf Proof.}\\
We use Blumenthal Theorem \ref{1.5}. $\lim_{n \to \infty} B_n = 0$ and $\lim_{n \to \infty} C_n =0$. Therefore $\sigma=0$ and $\tau =0$. Hence the property holds by using Remark \ref{1.6}.
$\square$

%**************************************************************************************************************************

%_____________________________________________section 3.3__________________________

%**************************************************************************************************************************
\subsection{Jacobi case (d)} 

The measure associated to $c_0$ is the Jacobi measure $(1-x)^\alpha (1+x)^\beta dx$ with $\alpha >-1$ and $\beta >-1$, and the one associated to $c_1$ is  $\dfrac{(1-x)^{\alpha+1} (1+x)^{\beta+1}}{\mid x-\xi \mid} dx+M \delta(\xi)$ with $\mid \xi \mid \ge 1$ and $M \ge 0$. Therefore the monic polynomials $P_n(x)$,  orthogonal with respect to $c_0$, are  the monic Jacobi polynomials $P_n^{(\alpha,\beta)}(x)$. Also here and for the same reasons as in the Laguerre case (c) we will denote by $c^{(\alpha,\beta)}$ the linear functional $c_0$. For more convenience $\mid x-\xi \mid$ will be written as $\varepsilon (x-\xi)$ with $\varepsilon=1$ if $\xi \le -1$, and $\varepsilon=-1$ if $\xi \ge 1$. Let us denote by $T_n(x,\xi,M)$, $n\ge0$, the monic polynomials orthogonal with respect to $c_1$. Moreover we will denote by $\tilde c_1$ the linear functional associated to the measure $\dfrac{(1-x)^{\alpha+1} (1+x)^{\beta+1}}{\varepsilon (x-\xi)} dx$ and therefore, we have $c^{(\alpha+1,\beta+1)}(.)=\tilde c_1(\varepsilon(x-\xi) .)$.\\
Relation (\ref{eq:1.1}) becomes
\begin{equation}\label{eq:J1.2}
T_n(x,\xi,M)= P_n^{(\alpha+1,\beta+1)}(x) - \sigma_n P_{n-1}^{(\alpha+1,\beta+1)}(x).
\end{equation}
Like in the Laguerre case (c) we have a property giving a relation between $T_{n+1}(t,\xi,M)$ and the reproducing kernel $N_n^{(\alpha+1,\beta+1)}(x,t)$.
%property
 \begin{property}\label{3.9}
  The monic polynomial $T_{n+1}(t,\xi,M)$, orthogonal with respect to $c_1$, satisfies the following relation
 \begin{equation}\label{eq:J1.3}
 T_{n+1}(t,\xi,M)= -\frac{k_n^{\alpha+1}}{c_1(P_n^{(\alpha+1,\beta+1)}(x))} \left(1-(t-\xi)c_1(N_n^{(\alpha+1,\beta+1)}(x,t))\right)
  \end{equation}
 where the linear functional $c_1$ acts on the variable $x$ and $N_n^{(\alpha+1,\beta+1)}(x,t)$ is the reproducing kernel 
  \begin{equation}\label{eq:J2.1}
 N_n^{(\alpha+1,\beta+1)}(x,t)= \sum_{j=0}^n \frac{P_j^{(\alpha+1,\beta+1)}(x)P_j^{(\alpha+1,\beta+1)}(t)}{k_j^{(\alpha+1,\beta+1)}}.
 \end{equation}
  \end{property}
 \noindent {\bf Proof.}\\
 $T_{n+1}(x,\xi,M)$ is written in the basis $\{1,\{\varepsilon(x-\xi)P_j^{(\alpha+1,\beta+1)}(x)\}_{j=0}^n\}$. Thus
 \begin{equation}\label{eq:J2.2}
 T_{n+1}(x,\xi,M)=T_{n+1}(\xi,\xi,M)+ \sum_{j=0}^n \theta_{j,n} \varepsilon(x-\xi)P_j^{(\alpha+1,\beta+1)}(x), \quad \theta_{n,n}=1.
 \end{equation}
 The sequel of the proof is the same as in Property \ref{2.12}. In the Jacobi case (d) we find 
  \begin{equation}\label{eq:J2.3}
  \theta_{i,n}=- \frac{T_{n+1}(\xi,\xi,M)}{k_i^{(\alpha+1,\beta+1)}}c_1(P_i^{(\alpha+1,\beta+1)}(x) ), \quad i=0,\ldots,n,
 \end{equation}
and since $\theta_{n,n}=1$ we have
 \begin{equation}\label{eq:J2.2a}
T_{n+1}(\xi,\xi,M)=- \frac{k_n^{(\alpha+1,\beta+1)}}{c_1(P_n^{(\alpha+1,\beta+1)}(x))}.
 \end{equation}
$\square$\\
The property giving $\hat k_n$ and $\sigma_n$ also has a similar proof as the one in Property \ref{2.16}.
%property
\begin{property}\label{3.12}
The square norm $\hat k_n$ of the polynomials $T_{n}(x,\xi,M)$ and $\sigma_n$ are given by
 \begin{eqnarray}
\hat k_0&=&c_1(1), \nonumber\\
\hat k_n&=& - k_{n-1}^{(\alpha+1,\beta+1)} \frac{c_1(P_n^{(\alpha+1,\beta+1)}(x))}{c_1(P_{n-1}^{(\alpha+1,\beta+1)}(x))}, \quad \forall nÊ\ge 1,\label{eq:J3.3}\\
\sigma_n &=& \frac{c_1(P_n^{(\alpha+1,\beta+1)}(x))}{c_1(P_{n-1}^{(\alpha+1,\beta+1)}(x))},  \quad \forall nÊ\ge 1. \label{eq:J3.4}
\end{eqnarray}
\end{property}
At last we have the following theorem concerning the Markov-Bernstein inequalities.
%theorem
\begin{theorem}\label{3.13}
The following Markov-Bernstein inequality holds:
$$
c_1((p^\prime)^2) \le \frac{1}{\mu_{1,n}} c_0(p^2), \quad \forall p \in \mathcal{P}_n
$$
where $\mu_{1,n}$ is the smallest zero of the polynomials $A_n(\lambda)$ satisfying the following three term recurrence relation
\begin{equation}\label{eq:J4.1}
A_n (\lambda) = (\lambda +B_n)A_{n-1} (\lambda)-C_n A_{n-2} (\lambda), \quad \forall n\ge 2.
\end{equation}
where 
\begin{eqnarray}
B_n&=& \frac{4(n-1)(n+\alpha)(n+\beta)}{n(2 n +\alpha+\beta-1)(2 n +\alpha+\beta)^2(2 n +\alpha+\beta+1)} \frac{c_1(P_{n-2}^{(\alpha+1,\beta+1)}(x))}{c_1(P_{n-1}^{(\alpha+1,\beta+1)}(x))} \nonumber \\
&&+ \frac{1}{(n-1)(n+\alpha+\beta)}
\frac{c_1(P_{n-1}^{(\alpha+1,\beta+1)}(x))}{c_1(P_{n-2}^{(\alpha+1,\beta+1)}(x))}, \label{eq:J4.2}\\
C_n&=& \frac{(n-2)(n+\alpha-1)(n+\beta-1)}{4(n-1)^2(n +\alpha+\beta)(2 n +\alpha+\beta-1)(2 n +\alpha+\beta-2)^2(2 n +\alpha+\beta-3)} \times \nonumber \\
&& \frac{c_1(P_{n-1}^{(\alpha+1,\beta+1)}(x)) c_1(P_{n-3}^{(\alpha+1,\beta+1)}(x))}{(c_1(P_{n-2}^{(\alpha+1,\beta+1)}(x)))^2}.\label{eq:J4.3}
\end{eqnarray}
$A_0(\lambda)=1$ and 
\begin{equation}\label{eq:J4.4}
A_1(\lambda)=\lambda- \frac{ k_1^{(\alpha+1,\beta+1)}}{c_1(1)} = \lambda-\frac{ 2^{\alpha+\beta+3} \Gamma(\alpha +2) \Gamma(\beta +2)}{(\alpha+\beta+2)_2 \Gamma(\alpha+\beta +3)c_1(1)}.
\end{equation}
\end{theorem}
\noindent {\bf Proof.}\\
By using (\ref{eq:J3.3}), (\ref{eq:J3.4}) and (\ref{eq:7.2}) we have
\begin{eqnarray*}
B_n&=& - \frac{k_{n}^{(\alpha,\beta)}}{n^2 k_{n-1}^{(1)}}- \frac{\sigma_{n-1}^2 }{({n-1})^2}\frac{ k_{n-1}^{(\alpha,\beta)}}{ k_{n-1}^{(1)}}\\
&=& \frac{4(n-1)(n+\alpha)(n+\beta)}{n(2 n +\alpha+\beta-1)(2 n +\alpha+\beta)^2(2 n +\alpha+\beta+1)} \frac{c_1(P_{n-2}^{(\alpha+1,\beta+1)}(x))}{c_1(P_{n-1}^{(\alpha+1,\beta+1)}(x))} \\
&&+ \frac{1}{(n-1)(n+\alpha+\beta)}
\frac{c_1(P_{n-1}^{(\alpha+1,\beta+1)}(x))}{c_1(P_{n-2}^{(\alpha+1,\beta+1)}(x))}, \\
C_n&=&  \frac{\sigma_{n-1}^2 (k_{n-1}^{(\alpha,\beta)})^2}{(n-1)^4  k_{n-1}^{(1)} k_{n-2}^{(1)}} \\
&=&  \frac{(n-2)(n+\alpha-1)(n+\beta-1)}{4(n-1)^2(n +\alpha+\beta)(2 n +\alpha+\beta-1)(2 n +\alpha+\beta-2)^2(2 n +\alpha+\beta-3)} \times  \\
&& \frac{c_1(P_{n-1}^{(\alpha+1,\beta+1)}(x)) c_1(P_{n-3}^{(\alpha+1,\beta+1)}(x))}{(c_1(P_{n-2}^{(\alpha+1,\beta+1)}(x)))^2}.
\end{eqnarray*}
$\square$\\
The Jacobi function of the second kind $Q_{n}^{(\alpha+1,\beta+1)}(\xi)$ can be used to replace the terms $c_1(P_{n}^{(\alpha+1,\beta+1)}(x))$ in $B_n$ and $C_n$ given in Theorem \ref{3.13}. $Q_{n}^{(\alpha,\beta)}(\xi)$ is defined as follows (see Szeg\"o \cite{S} Relations (4.61.4) and (4.61.5) page 74) for $\xi \in \mathbb{C} \setminus [-1,1]$:
\begin{eqnarray}
Q_{n}^{(\alpha,\beta)}(\xi)&=&\frac{1}{2}(\xi-1)^{-\alpha}(\xi+1)^{-\beta}\int_{-1}^1(1-x)^\alpha (1+x)^\beta \frac{\tilde P_{n}^{(\alpha,\beta)}(x)}{\xi-x} dx \label{eq:JJ.1}\\
&=&2^{n+\alpha +\beta} \frac{\Gamma(n+\alpha +1)\Gamma(n+\beta +1)}{\Gamma(2n+\alpha + \beta+2)} \frac{F(n+\alpha +1,n+1;2n+\alpha + \beta+2;\frac{2}{1-\xi})}{(\xi-1)^{n+\alpha +1}(\xi+1)^\beta} \label{eq:JJ.2}
\end{eqnarray}
where $\tilde P_{n}^{(\alpha,\beta)}(x)$ is the non-monic Jacobi polynomial ($\tilde P_{n}^{(\alpha,\beta)}(x)=\frac{(n+\alpha + \beta+1)_n}{2^n n!} P_{n}^{(\alpha,\beta)}(x)$) and $F$ is the hypergeometric function $F(a,b;c;z)=\sum_{\nu=0}^\infty \frac{(a)_\nu (b)_\nu}{(c)_\nu \nu!} z^\nu$.\\
Now the following property can be proved.
%property
\begin{property}\label{JJ1}
$$
\lim_{n \to \infty} \mu_{1,n}=0.
$$
\end{property}
\noindent {\bf Proof.}\\
First we prove this result when $\mid \xi \mid >1$. Then we will use the asymptotic behavior of $Q_{n}^{(\alpha+1,\beta+1)}(\xi)$ and of $P_{n}^{(\alpha+1,\beta+1)}(\xi)$ when $n$ tends to infinity. We have (see Szeg\"o \cite{S} Relation (8.71.19)  page 225)
$$
(\xi-1)^{\alpha +1}(\xi+1)^{\beta +1}Q_{n}^{(\alpha+1,\beta+1)}(\xi)\simeq\frac{1}{\sqrt {n}}(\xi-\sqrt {\xi^2-1})^{n+1} \Phi(\xi)
$$
where $\Phi$ is independent of $n$, and (see Szeg\"o \cite{S} Theorem 8.21.7  page 196)
$$
\tilde P_{n}^{(\alpha+1,\beta+1)}(\xi)\simeq \frac{(\sqrt {\xi-1}+\sqrt {\xi+1})^{\alpha+\beta +2}}{(\xi-1)^{(\alpha +1)/2}(\xi+1)^{(\beta +1)/2}}\frac{(\xi+\sqrt {\xi^2-1})^{n+1/2}}{\sqrt {2 \pi n}\sqrt {\xi^2-1}}.
$$
Therefore
\begin{eqnarray*}
c_1(P_{n}^{(\alpha+1,\beta+1)}(x))&=& \frac{2^n n!}{(n+\alpha + \beta+3)_n}(\tilde c_1(\tilde P_{n}^{(\alpha+1,\beta+1)}(x))+ M \tilde P_{n}^{(\alpha+1,\beta+1)}(\xi))\\
&=& \frac{2^n n!}{(n+\alpha + \beta+3)_n}(2 \varepsilon (\xi-1)^{\alpha +1}(\xi+1)^{\beta +1}Q_{n}^{(\alpha+1,\beta+1)}(\xi) \\
&&+ M \tilde P_{n}^{(\alpha+1,\beta+1)}(\xi)).
\end{eqnarray*}
For $\xi$ fixed, $\mid \xi \mid >1$, we have $\mid \xi-\sqrt {\xi^2-1} \mid <1$ and $\lim_{n \to \infty} \mid \xi-\sqrt {\xi^2-1} \mid =0$. Therefore, if $M>0$,
$$
c_1(P_{n}^{(\alpha+1,\beta+1)}(x)) \simeq \frac{2^n n! (\xi+\sqrt {\xi^2-1})^{n+1/2}}{(n+\alpha + \beta+3)_n \sqrt {n}}\left(\frac{M(\sqrt {\xi-1}+\sqrt {\xi+1})^{\alpha+\beta +2}}{(\xi-1)^{\alpha /2}(\xi+1)^{\beta /2}\sqrt {2 \pi }(\xi^2-1)} \right)
$$
and if $M=0$, 
$$
c_1(P_{n}^{(\alpha+1,\beta+1)}(x)) \simeq \frac{\varepsilon 2^{n+1} n!  (\xi-\sqrt {\xi^2-1})^{n+1}}{(n+\alpha + \beta+3)_n \sqrt {n}} \Phi(\xi).
$$
Hence, if $M>0$, we have
$$
\frac{c_1(P_{n-1}^{(\alpha+1,\beta+1)}(x))}{c_1(P_{n-2}^{(\alpha+1,\beta+1)}(x))} \simeq \frac{2\sqrt{(n-1)(n-2)}(n+\alpha + \beta+1)}{(2n+\alpha + \beta-1)_2}(\xi+\sqrt {\xi^2-1}).
$$
If $M=0$, we have
$$
\frac{c_1(P_{n-1}^{(\alpha+1,\beta+1)}(x))}{c_1(P_{n-2}^{(\alpha+1,\beta+1)}(x))} \simeq \frac{2\sqrt{(n-1)(n-2)}(n+\alpha + \beta+1)}{(2n+\alpha + \beta-1)_2}(\xi-\sqrt {\xi^2-1}).
$$
In the two cases we have $\lim_{n \to \infty} B_n = \lim_{n \to \infty} C_n=0$. Therefore $\sigma=\tau=0$ and the property holds by using Remark \ref{1.6}.
\medskip

When  $\mid \xi \mid =1$, $\tilde c_1=c^{(\alpha,\beta+1)}$ if $\xi=1$ and $\tilde c_1=c^{(\alpha+1,\beta)}$ if $\xi=-1$.\\
If $\xi=1$, by using recursively (\ref{eq:L9}) for $n,n-1,\ldots,1$ we obtain 
$$
\tilde c_1(P_{n}^{(\alpha+1,\beta+1)}(x))= c^{(\alpha,\beta+1)}(P_{n}^{(\alpha+1,\beta+1)}(x))=\frac{2^n n! (\beta+2)_n}{(\alpha+\beta+3)_{2n}} k_0^{(\alpha,\beta+1)}.
$$
Then
\begin{eqnarray*}
c_1(P_{n}^{(\alpha+1,\beta+1)}(x)) &=& \frac{2^{n+\alpha+\beta+2}n! (\beta+2)_n \Gamma(\alpha+1) \Gamma(\beta+2)}{(\alpha+\beta+3)_{2n}\Gamma(\alpha+\beta+3)}+MP_{n}^{(\alpha+1,\beta+1)}(1)\\
&=& \frac{2^{n} (\alpha+2)_n}{(n+\alpha+\beta+3)_{n}}(M+2^{\alpha+\beta+2}\frac{n!}{(\alpha+2)_n} \frac{(\beta+2)_n}{(\alpha+\beta+3)_{n}} \frac{ \Gamma(\alpha+1) \Gamma(\beta+2)}{\Gamma(\alpha+\beta+3)}).
\end{eqnarray*}
Since $\alpha+1>0$ and $\beta+2>1$ the following quantities tend to 0 when $n$ tends to $\infty$:
\begin{eqnarray*}
\frac{n!}{(\alpha+2)_n}&=&  \prod_{j=1}^n \frac{1}{1+\frac{\alpha+1}{j}},\\
 \frac{(\beta+2)_n}{(\alpha+\beta+3)_{n}}&=& \prod_{j=1}^n \frac{1}{1+\frac{\alpha+1}{\beta +1+j}}.
\end{eqnarray*}
Therefore, if $M>0$,
$$
c_1(P_{n}^{(\alpha+1,\beta+1)}(x)) \simeq M  \frac{2^{n} (\alpha+2)_n}{(n+\alpha+\beta+3)_{n}}
$$
and if $M=0$, 
$$
c_1(P_{n}^{(\alpha+1,\beta+1)}(x)) =  \frac{2^{n+\alpha+\beta+2}n! (\beta+2)_n \Gamma(\alpha+1) \Gamma(\beta+2)}{(\alpha+\beta+3)_{2n}\Gamma(\alpha+\beta+3)}.
$$
Hence, if $M>0$, we have
$$
\frac{c_1(P_{n-1}^{(\alpha+1,\beta+1)}(x))}{c_1(P_{n-2}^{(\alpha+1,\beta+1)}(x))} \simeq \frac{2(n+\alpha)(n+\alpha + \beta+1)}{(2n+\alpha + \beta-1)_2}.
$$
If $M=0$, we have
$$
\frac{c_1(P_{n-1}^{(\alpha+1,\beta+1)}(x))}{c_1(P_{n-2}^{(\alpha+1,\beta+1)}(x))} \simeq \frac{2(n+ \beta)}{(2n+\alpha + \beta-1)_2}.
$$
Again, in the two cases we have $\lim_{n \to \infty} B_n = \lim_{n \to \infty} C_n=0$. Therefore $\sigma=\tau=0$ and the property holds by using Remark \ref{1.6}.\\
If $\xi=-1$, by using (\ref{eq:L8}) all the sequel is analogous.
$\square$
%**************************************************************************************************************************

%_____________________________________________section 4_________________________

%**************************************************************************************************************************
\section{Conclusion} 

The seven kinds of three term recurrence relation have been given. If $A_n(\lambda)$ is a polynomial obtained for one of them, the inverse of the square root of its smallest zero is the Markov-Bernstein constant linked to the corresponding coherent pair of measures.\\
Five kinds of three term recurrence relation are totally explicit. The last two relations are provided by means of integrals: $c_1(L_{j}^{(\alpha+1)}(x))$ in the Laguerre case (c) and $c_1(P_{j}^{(\alpha+1,\beta+1)}(x))$ in the Jacobi case (d) which can possibly be replaced by functions of the second kind.\\
In \cite{DEH2} and \cite{DRA} classical numerical methods, in particular the Newton method and the $qd$ algorithm, were used to obtain explicit upper and lower bounds of Markov-Bernstein constants. Here these methods have only been used successfully in the Laguerre case (b). Even in the case of the very simple three term recurrence relation (\ref{eq:10.1}) obtained in the Laguerre case (a) these methods do not provide an explicit result. But we have proved in the Laguerre case (a), the Laguerre case (b) when $\xi=0$ and all the Jacobi cases that $\lim_{n \to \infty} \mu_{1,n}=0$.
Nevertheless these three term recurrence relations are the prerequisite in any advanced study of Markov-Bernstein constants linked to coherent pairs of measures. In particular the numerical value of $\mu_{1,n}$ can be computed from this three term recurrence relation in order to use it in some applications.

%\bibliographystyle{amsplain}
%\bibliography{R\'ef\'erences}

\end{document}